\documentclass[11pt]{article}

\usepackage{amsmath,amsthm,amssymb,enumerate,dsfont}
\usepackage{algorithmicx,algorithm,algpseudocode}
\usepackage{float}
\usepackage{titlesec}
\usepackage{scalerel}
\usepackage{epstopdf}
\usepackage{graphicx}
\usepackage{xspace}
\usepackage{color}
\usepackage{hyperref}
\usepackage[title]{appendix}
\usepackage[left=2.5cm,right=2.5cm,vmargin=2.5cm]{geometry}
\usepackage[numbers,sort&compress]{natbib}
\usepackage{caption}
\usepackage{subcaption}
\usepackage[labelformat=simple]{subcaption}
\usepackage{bm}
\usepackage{cases}
\usepackage{multirow}

\newcommand{\MATLAB}{\textsc{Matlab}\xspace}

\newcommand{\G}{\mathcal{G}}
\newcommand{\R}{\mathbb{R}}
\newcommand{\F}{\mathcal{F}}
\newcommand{\SO}{\mathcal{SO}}
\newcommand{\OO}{\mathcal{O}}
\newcommand{\tr}{\mathrm{tr}}
\newcommand{\ip}[1]{\left<{#1}\right>}

\newcommand{\hU}{\widehat{\bm{U}}}
\newcommand{\RK}{\left(\mathbf{1}_{1\times K}\otimes\bm R\right)}
\newcommand{\setbar}{\;\middle|\;}

\usepackage[colorinlistoftodos,prependcaption,textsize=tiny]{todonotes}

\newcommand{\nm}[1]{\left\|{#1}\right\|}

\newcommand{\ie}{{\it i.e.}}

\DeclareMathOperator*{\argmin}{argmin}
\DeclareMathOperator*{\argmax}{argmax}

\newtheorem{lem}{Lemma}
\newtheorem{thm}{Theorem}

\newtheorem{defi}{Definition}
\newtheorem{prop}{Proposition}

\newtheorem{rmk}{Remark}

\begin{document}
\title{\textbf{Non-Convex Joint Community Detection and Group Synchronization via 
Generalized Power Method}}
\author{Sijin Chen\thanks{Department of Computer Science and Engineering, the Chinese University of Hong Kong (CUHK). E--mail: {\tt sjchen0@cse.cuhk.edu.hk}}
\and Xiwei Cheng\thanks{Department of Computer Science and Engineering, CUHK.  E--mail: {\tt xwcheng0@cse.cuhk.edu.hk}}
\and Anthony Man-Cho So\thanks{Department of Systems Engineering and Engineering Management, CUHK. E--mail: {\tt manchoso@se.cuhk.edu.hk}}}
\date{ }
\maketitle
\begin{abstract}
This paper proposes a \textit{Generalized Power Method} (GPM) to tackle the problem of community detection and group synchronization simultaneously in a direct non-convex manner. Under the \textit{stochastic group block model} (SGBM), theoretical analysis indicates that the algorithm is able to exactly recover the ground truth in $O(n\log^2 n)$ time, sharply outperforming the benchmark method of semidefinite programming (SDP) in $O(n^{3.5})$ time. Moreover, a lower bound of parameters is given as a necessary condition for exact recovery of GPM. The new bound breaches the information-theoretic threshold for pure community detection under the stochastic block model (SBM), thus demonstrating the superiority of our simultaneous optimization algorithm over the trivial two-stage method which performs the two tasks in succession. We also conduct numerical experiments on GPM and SDP to evidence and complement our theoretical analysis.
\end{abstract}

\bigskip

\noindent {\bf Keywords:} Community detection, group synchronization, Generalized Power Method, iterative algorithms, stochastic block model, orthogonal group, rotational group
\bigskip

\section{Introduction}
Community detection methods typically make use of the edge connectedness information of a given observation network to generate an inference on the underlying clustering of the agents \cite{AbbeCommunityDetection2018,yun2014accurate,amini2016semidefinite,gao2015achieving}. However, if additional information of the agents besides the edge connectedness is available, chances are that the recovery results of the clustering can outperform the pure community detection methods, even its information-theoretic limit \cite{Abbe2015CommunityDI}, by carefully exploiting these extra information \cite{weng2016community,binkiewicz2017covariate,zhang2016community}. Based on this fact, one may proceed to wonder whether a simultaneous recovery of both the underlying clustering and the additional agent structures can be achieved efficiently, or, more efficiently than the trivial two-stage method that classifies all the agents first and then recovers the structures according to the classification results. In this paper, we study this question by considering the scenarios where additional information comes from pairwise relative measurements in a group $\G$, which usually appears in group synchronization problems \cite{arie2012global, boumal2016nonconvex,bandeira2016low,bandeira2018random,liu2020unified}. This assumption is also motivated by a number of practical problems, for example the 2D class averaging in cryo-electron microscopy single particle reconstruction \cite{frank2006three,amit2011viewing,zhao2014rotationally}, where 2D images rotated in different directions are required to be classified and aligned at the same time before further processes to denoise the observation. 

\begin{figure}[h!]
    \centering
    \includegraphics[width=0.8\textwidth]{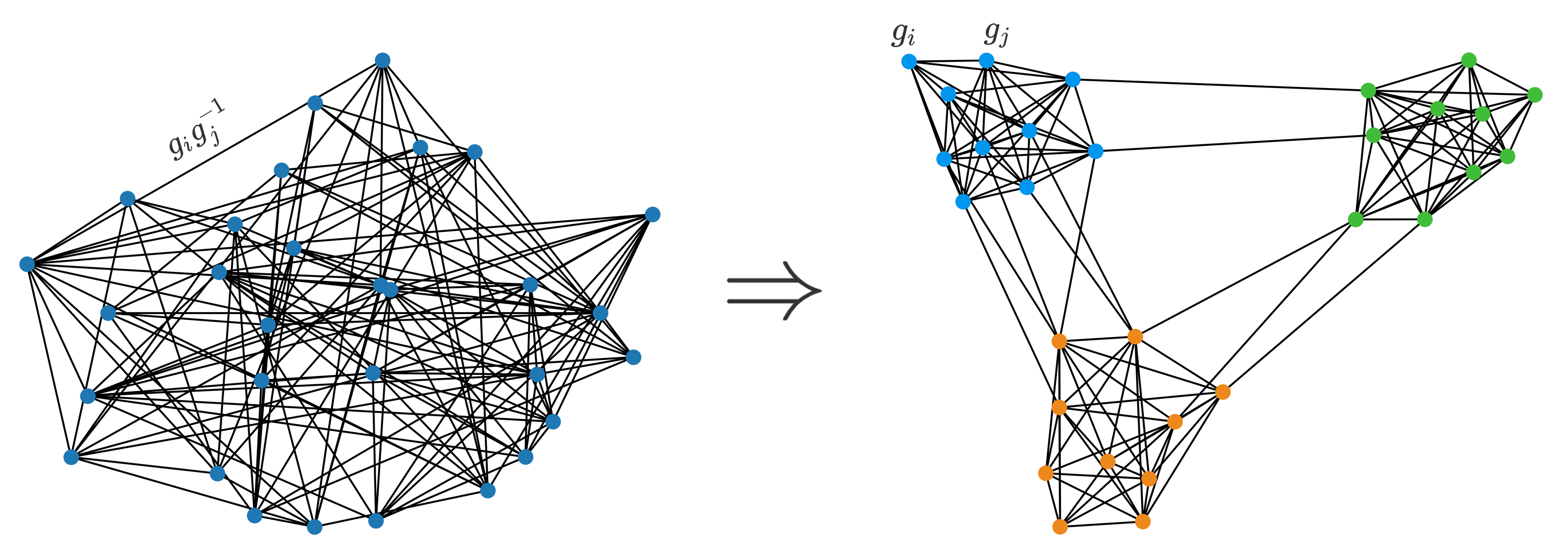}
    \caption{Illustration of the joint optimization problem of community detection and group synchronization. We observe a network of $n=30$ agents falling in $K=3$ equal-sized communities, with relative measurement in $\G$ between any connected pair of nodes. The target is to recover both the underlying cluster and the group element of each agent simultaneously.}
    \label{fig:demo}
\end{figure}

Recently, in analogy to the celebrated stochastic block model (SBM) \cite{Abbe2015CommunityDI} for community detection, \cite{fan2021joint} defined a model and its joint optimization problem of community detection and synchronization, which we inherit with necessary adjustments and generalizations. Assume that there are $n$ agents in a network partitioned into $K$ communities, and each agent $i$ corresponds to a group element $g_i\in\G$. With probability $p$, we obtain the relative group measurement $g_ig_j^{-1}\in\G$ for two agents $i,j$ belonging to the same cluster; with probability $q$, we obtain an observation noise $g$ uniformly sampled from $\G$ for two agents $i,j$ belonging to different clusters -- a mechanism named the outlier noise model \cite{Singer2011AngularSB,Fan2019multifrequency}. Given these observations, we aim to recover both the underlying clustering and the corresponding group elements of all the agents. Figure \ref{fig:demo} illustrates the problem settings and the target of recovery.

This problem can be formulated as a maximum likelihood estimation (MLE) problem, yet it is still non-convex and computationally challenging to solve. In \cite{fan2021joint}, the authors proposed a semidefinite programming (SDP) method to jointly perform the clustering and synchronization task. This joint algorithm first relaxes the original problem into a convex one that involves both the clustering variables and the rotational group variables, and then use SDP to solve the convex problem. By exploiting the mutual rotational measurements to recover the underlying community and rotation of each agent simultaneously, it enjoys both theoretical and experimental advantages over the state-of-the-art pure community detection methods. However, applying SDP to practice is yet time-consuming when the problem scale gets large, since the standard interior point method for SDP typically takes $O(n^{3.5})$ time to find an optimal solution. The theory developed in \cite{fan2021joint} also leaves several open questions. For example, the theory only takes care of the $\G=\SO(d)$ situation, and rigorous analysis is only provided for two-cluster cases. These issues hinder the reliability of SDP when applied to general situations.

\subsection{Contributions}
Different from the methodology of convex relaxation adopted by the authors of \cite{fan2021joint}, this paper proposes an iterative \textit{Generalized Power Method} (GPM) to directly tackle the non-convex optimization problem of simultaneous community detection with group synchronization, and provides a theoretical guarantee for its linear convergence to the optimal solution under certain conditions. Our contributions are threefold.

\begin{itemize}
    \item Most significantly, our algorithm sharply reduces the time complexity to $O(n\log^2 n)$ from $O(n^{3.5})$ of SDP without any compensation on the lower bound for parameters. The iterative algorithm is also structurally simple and practically convenient to implement.
    
    \item Moreover, in this paper, theoretical guarantees are provided for both rotational ($\G=\SO(d)$) and orthogonal ($\G=\OO(d)$) synchronization, and the number of clusters can be finitely many, \ie, $\Theta(1)$. These generalize the scenario of rotational synchronization with only 2 clusters, to which the theoretical results of \cite{fan2021joint} apply. 
    \item We also remark that the conditions for linear convergence in this paper are able to break the information-theoretic lower bound $\sqrt{\alpha}-\sqrt{\beta}>\sqrt{K}$ for pure community detection \cite{Abbe2015CommunityDI}, thus demonstrating the superiority of our joint method over the naive two-stage approaches.
\end{itemize}

\subsection{Organization}
The rest of this paper is organized as follows. Section \ref{sec:preliminaries} briefly introduces the preliminaries about community detection and group synchronization, which are important to the presentation and analysis of our joint problem. Then, Section \ref{sec:formulation} presents a formal definition of our probabilistic model and formulate the non-convex optimization problem of simultaneous community detection with group synchronization. Our non-convex methodology is elaborated in Section \ref{sec:methodology}, followed by a detailed statement of our main theoretical results in Section \ref{sec:main-result}. Their proofs are developed in Section \ref{sec:main-result}--\ref{sec:spec-init}. Section \ref{sec:numerics} presents the numerical results from computer simulations. Finally, we make some discussions and conclude this paper in Section \ref{sec:conc}.

\subsection{Notations}
Throughout this paper, all matrices, vectors, and scalars are assumed to be real. We follow the convention of using boldface letters to represent matrices and vectors, and reserve the normal letters for scalars. We denote $\bm X^\top$ the transpose of a matrix $\bm X$. If the parameter $d$ (to be defined later) divides the size of matrix $\bm X$, we refer $\bm X_{ij}$ to the $(i,j)$-th block of $\bm X$ of size $d\times d$ unless otherwise specified, and refer $\bm X_{i\times}$ to the $i$-th block row of $\bm X$. 

For any matrix $\bm X\in\R^{n\times n}$ and any positive integer $k\leq n$, $\lambda_k(\bm X)$ is the $k$-th largest singular value of $\bm X$ so that $\{\sigma_k(\bm X)\}$ forms a non-ascending sequence. Similarly, $\lambda_k(\bm X)$ denotes its $k$-th largest (real) eigenvalue, given that $\bm X$ is symmetric. $\nm{\bm X}:=\sigma_1(\bm X)$, $\nm{\bm X}_F:=\tr\left(\bm X^\top\bm X\right)$, and $\nm{\bm X}_*:=\sum\sigma_k(\bm X)$ respectively represents the operator norm, Frobenius norm, and nuclear norm of $\bm X$. We denote $\ip{\bm X, \bm Y}:=\tr\left(\bm X^\top\bm Y\right)$ the (Frobenius) inner product of two matrices $\bm X$ and $\bm Y$. 

Consider the spectral decomposition $\bm Q\bm \Lambda\bm Q^\top$ for any given symmetric $\bm X\in\R^{n\times n}$, where the diagonal elements in $\bm\Lambda$ are non-ascending. Based on \MATLAB notation, we define a mapping $\bm U=\mathrm{eigs}_{[k:l]}(\bm X)$, such that $\bm U$ is the submatrix containing the $k$-th to $l$-th column of $\bm Q$, where $1\leq k\leq l\leq n$. When $k=1$, we write for simplicity that $\bm U=\mathrm{eigs}_{l}(\bm X)$.

For any integer $n,m\geq 1$, the $n\times n$ identity matrix is denoted by $\bm I_{n}$, the $n\times m$ all-one matrix is denoted by $\mathbf{1}_{n\times m}$, and the $n\times m$ all-zero matrix is denoted by $\mathbf{0}_{n\times m}$. We use $\otimes$ and $\odot$ to represent the operator of Kronecker product and Hadamard (elementwise) product of two matrices with conforming shapes.

For two nonnegative functions $f(n)$ and $g(n)$, we say that $f(n)=O(g(n))$ if there exists some $C$ and $N$ such that $f(n)\leq Cg(n)$ for all $n\geq N$, and $f(n)=\Omega(g(n))$ if there exists some $C$ and $N$ such that $f(n)\geq Cg(n)$ for all $n\geq N$. We say that $f(n)=\Theta(g(n))$ if $f(n)=O(g(n))$ and $f(n)=\Omega(g(n))$, and $c$ is a constant if $c=\Theta(1)$. Also, $f(n)=o(g(n))$ if $\lim_{n\to\infty}\frac{f(n)}{g(n)}=0$. For any integer $n\geq 1$ we define $[n]:=\{1,2,...,n\}$.

\section{Preliminaries}
This section introduces some important ingredients of community detection and group synchronization problems that serve as preliminaries for the development of the remaining contents.

\subsubsection*{Community detection}
The main problem to be studied in this paper has a strong correlation to \textit{community detection} problems under the symmetric stochastic block model (SBM) \cite{Abbe2015CommunityDI} with parameters $(n,p,q,K)$. Assume that $n$ agents in a network fall in $K$ underlying communities of equal size $m=n/K$ (balanced clustering). Then, $\mathrm{SBM}(n,p,q,K)$ generates a random undirected graph $G$ such that each two distinct nodes are connected by an edge with probability $p$ if they belong to the same cluster, and with probability $q$ otherwise. We assume without loss of generality that a node is always connected to itself in $\mathrm{SBM}(n,p,q)$. 

Clustering functions and clustering matrices \cite{wang2021optimal} are defined to formally represent the community structure of the nodes. We denote $C:[n]\to[K]$ the clustering function that maps node $i$ to cluster $C(i)$ where it belongs. Conversely, we define $\mathcal{I}_j:=\left\{i\in[n]\;\middle|\;C(i)=j\right\}$, that is, the pre-image of $C$ at cluster $j$. A matrix $\bm H\in\{0,1\}^{n\times K}$ is said to be a clustering matrix, if there exists a clustering function $C(i)$ such that $\bm H_{ij}=1$ if and only if $C(i)=j$. This establishes a one-one correspondence between $C$ and $\bm H$. 

Further, one can show that $\bm H\in\{0,1\}^{n\times K}$ is a clustering matrix if and only if $\bm H\mathbf{1}_{K\times 1}=\mathbf{1}_{n\times 1}$, and $\mathbf{1}_{1\times n}\bm H=m\mathbf{1}_{1\times K}$. It is then natural to define a set $\mathcal{H}$ of all $n\times K$ clustering matrices, such that 
$$\mathcal{H}:=\left\{\bm H\in\{0,1\}^{n\times K}\;\middle|\;\bm H\mathbf{1}_{K\times 1}=\mathbf{1}_{n\times 1}, \mathbf{1}_{1\times n}\bm H=m\mathbf{1}_{1\times K}\right\}.$$
The fact that any permutation of the numbering of the clusters does not matter gives an equivalence relation on $\mathcal{H}$. For any $\bm H\in\mathcal{H}$, the equivalent class of $\bm H$ is defined to be $\left\{\bm H\bm P\;\middle|\;\bm P\in\mathcal{S}_K\right\}$, where $\mathcal{S}_K$ is the $K$-dimensional permutation group. Therefore, given a ground truth $\bm H^*$, the estimation error \cite{wang2021optimal} of $\bm H\in\mathcal{H}$ is defined as 
$$\epsilon(\bm H)=\min_{\bm P\in\mathcal{S}_K}\nm{\bm H-\bm H\bm P}_F.$$

For arbitrary matrix $\bm M\in\R^{n\times K}$, we say
\begin{equation}\label{pr:proj-H}
    \Pi_\mathcal{H}(\bm M):=\argmin_{\bm H\in\mathcal{H}}\nm{\bm M-\bm H}_F=\argmax_{\bm H\in\mathcal{H}}\ip{\bm H,\bm M}
\end{equation}
is the projection of $\bm M$ onto $\bm H$, where the second equality is due to the fact that $\nm{\bm H}_F=\sqrt{n}$ for all $\bm H\in\mathcal{H}$. As is pointed out by \cite{wang2021optimal}, Problem \ref{pr:proj-H} is equivalent to a \textit{minimum-cost assignment problem} (MCAP) that can be tackled efficiently by existing algorithms, such as \cite{Tokuyama1995}. The following proposition will be useful in the later analysis on the time complexity of our algorithm.
\begin{prop}[Proposition 1, \cite{wang2021optimal}]\label{prop:runtime-MCAP}
    Problem \ref{pr:proj-H} can be solved in $O(K^2n\log n)$ time.
\end{prop}

When the parameters $p$ and $q$ are located in the logarithmic sparsity region of $\mathrm{SBM}(n,p,q,K)$, \ie, $p=\frac{\alpha\log n}{n}$ and $\frac{\beta\log n}{n}$ where $\alpha,\beta=O(1)$, \cite{Abbe2015CommunityDI} derived that one can recover the underlying clustering if and only if $\sqrt{\alpha}-\sqrt{\beta}>\sqrt{K}$ under SBM. This is the information-theoretic limit for the pure community detection problems.

\subsubsection*{Group synchronization}
The main problem interacts with another non-convex optimization problem named \textit{group synchronization}, and our focus will be on the synchronization of orthogonal and rotational groups. We denote $\OO(d)$ the $d$-dimensional orthogonal group over $\R$, such that 
$$\OO(d):=\left\{\bm Q\in\R^{d\times d}\;\middle|\;\bm Q\bm Q^\top=\bm Q^\top\bm Q=\bm I_d\right\}$$
with the usual matrix multiplication as the group operation. The $d$-dimensional rotational group $\SO(d)$, or special orthogonal group, is the kernel of the group homomorphism $\det(\cdot)$ on $\OO(d)$. In other words,
$$\SO(d):=\left\{\bm Q\in\R^{d\times d}\;\middle|\;\bm Q\bm Q^\top=\bm Q^\top\bm Q=\bm I_d,\;\det(\bm Q)=1\right\}.$$
In the typical formulation of group synchronization \cite{boumal2016nonconvex,liu2020unified}, there are $n$ agents in a measurement network and the $i$-th agent corresponds to a group element $g_i\in\G$, where we specially consider $\G\in\{\OO(d),\SO(d)\}$ in this paper. For group $\G$, we define a set of block matrices $\G^n\subset\R^{nd\times d}$ as follows:
\begin{equation*}
    \G^n:=\left\{\bm X\in\R^{nd\times d}\;\middle|\;\bm X=\begin{pmatrix}g_1\\g_2\\\vdots\\g_n\end{pmatrix},g_i\in\G\right\}.
\end{equation*}
Moreover, a block-diagonalization operator $\mathrm{bdiag}(\cdot)$ is defined, such that
\begin{equation*}
    \mathrm{bdiag}(\bm X):=\begin{pmatrix}
    g_1&&&\\
    &g_2&&\\
    &&\ddots&\\
    &&&g_n
    \end{pmatrix},\;\forall\bm X=\begin{pmatrix}
    g_1\\g_2\\\vdots\\g_n
    \end{pmatrix}\in\mathcal{G}^n,
\end{equation*}
and we also denote $\mathrm{bdiag}(\G^n):=\left\{\mathrm{bdiag}(\bm X)\setbar\bm X\in\G^n\right\}$.

The design and analysis of our algorithm heavily relies on the projection onto $\OO(d)$, which is also prevalent in various works on group synchronization problems \cite{arie2012global,ling2020nearoptimal,liu2020unified}. For arbitrary matrix $\bm X\in\R^{d\times d}$, we define 
\begin{equation}\label{pr:proj-Od}
    \Pi_{\OO(d)}(\bm X):=\argmin_{\bm Q\in\OO(d)}\nm{\bm X-\bm Q}_F=\argmax_{\bm Q\in\OO(d)}\ip{\bm X,\bm Q}
\end{equation}
the projection of $\bm X$ onto $\OO(d)$, where the second equality is due to the fact that $\nm{\bm Q}_F=\sqrt{d}$ for all $\bm Q\in\OO(d)$. Problem \ref{pr:proj-Od} can be categorized into the \textit{orthogonal Procrustes problem} \cite{GowerProcrustes2004} that has a closed-form solution as follows: 
\begin{prop}\label{prop:proj-Od}
    For any given $\bm X\in\R^{d\times d}$, suppose that $\bm U\bm \Sigma\bm V^\top$ is the SVD of $\bm X$. Then, $\Pi_{\OO(d)}(\bm X)=\bm U\bm V^\top$, and $\max_{\bm Q\in\OO(d)}\ip{\bm X, \bm Q}=\tr(\bm\Sigma)=\nm{\bm X}_*$.
\end{prop}\label{sec:preliminaries}

\section{Problem formulation}\label{sec:formulation}
In this section, our probabilistic model is defined and the concerning non-convex optimization problem is formulated as a maximum likelihood estimation (MLE). We then properly reformulate this problem to fit for our methodology and define a metric of estimation error.

\subsection{The probabilistic model and the joint optimization problem}

We now formally state our joint model named \textit{stochastic group block model} (SGBM) with parameters $(n,p,q,K,d,\G)$, and formulate the joint optimization problem as a maximum likelihood estimation (MLE). Assume that $n$ agents in a network fall in $K$ underlying communities of equal size $m=n/K$, and each agent $i$ corresponds to a group element $\bm R^*_i\in\G$, where $\G\in\{\OO(d),\SO(d)\}$. All the group elements constitute a block matrix $\bm R^*\in\G^n$. Denote $C^*$ the clustering function, $\mathcal{I}^*_k$ the pre-image of $C^*$ at cluster $k$, and $\bm H^*\in\mathcal{H}$ the clustering matrix. $\mathrm{SGBM}(n,p,q,K,d,\G)$ firstly generates a random undirected graph $G=(V,E)$ under $\mathrm{SBM}(n,p,q,K)$. Then, it generates an observation matrix $\bm A\in\R^{nd\times nd}$, where its $(i,j)$-th block
\begin{equation}\label{pr:setting}
\begin{split}
\bm A_{ij}=\begin{cases}
\bm R^*_i\bm R^{*\top}_j,&\text{if }\{i,j\}\in E\text{ and }C^*(i)=C^*(j),\\
\bm R_{ij}\sim\mathrm{\mathrm{Unif}(\G)},&\text{if }\{i,j\}\in E,i<j,\text{ and }C^*(i)\neq C^*(j),\\
\bm R_{ji}^\top ,&\text{if }\{i,j\}\in E,i>j,\text{ and }C^*(i)\neq C^*(j),\\
\bm 0,&\text{otherwise,}
\end{cases}
\end{split}
\end{equation}
where $\mathrm{Unif}(\G)$ is the uniform distribution over $\G$ with respect to the Haar measure. As an equivalent statement, $\mathrm{SGBM}(n,p,q,K,d,\G)$ generates a realization $\bm A$ of the random matrix $\bm A^\text{r}$, where its $(i,j)$-th block
\begin{equation}\label{pr:setting-equivalent}
\begin{split}
\bm A^\text{r}_{ij}=\begin{cases}
w_{ij}\bm R^*_i\bm R^{*\top}_j:w_{ij}\sim\mathrm{Bern}(p),&\text{if }i<j\text{ and }C^*(i)=C^*(j),\\
u_{ij}\bm R_{ij}:u_{ij}\sim\mathrm{Bern}(q),\bm R_{ij}\sim\mathrm{\mathrm{Unif}(\G)},&\text{if }i<j\text{ and }C^*(i)\neq C^*(j),\\
\bm A_{ji}^\top,&\text{if }i>j,\\
\bm I_d,&\text{otherwise.}
\end{cases}
\end{split}
\end{equation}

Given the observation matrix $\bm A$, \cite{fan2021joint} introduced the MLE problem for both the community structure $\bm H\in\mathcal{H}$ and the group elements $\bm R\in\G^n$, referred to as \textsf{Joint-}$\G$ in this paper:
\begin{equation}\label{pr:orig}
\max_{\substack{\bm R\in\G^n\\\bm H\in\mathcal{H}}}\sum_{i,j\in \mathcal{I}_k,\forall k}\left<\bm A_{ij}, \bm R_i\bm R_j^\top \right>.\tag{\textsf{Joint-}$\G$}
\end{equation}

\subsection{Reformulation of the problem}
Problem \ref{pr:orig} is computationally intractable without proper reformulations. Defining $\bm M\in\R^{nd\times nd}$ such that its $(i,j)$-block
\begin{equation}\label{eq:def-M}
\bm M_{ij}=
\begin{cases}
    \bm R_i\bm R_j^\top, &\text{if }C(i)=C(j),\\
    \mathbf{0}, &\text{otherwise},
\end{cases}
\end{equation}
we can convert Problem \ref{pr:orig} to
\begin{equation}\label{pr:SDP}
    \max_{\substack{\bm R\in\G^n\\\bm H\in\mathcal{H}}}\left<\bm A, \bm M\right>.
\end{equation}
Since 
\begin{equation*}
\bm H_{i\times}\bm H_{j\times}^\top=
\begin{cases}
    1, &\text{if }C(i)=C(j),\\
    0, &\text{otherwise},
\end{cases}
\end{equation*}
we have $\bm M_{ij}=\left(\bm H_{i\times}\bm H_{j\times}^\top\right)\bm R_i\bm R_j^\top$. Then, direct calculations can verify the following proposition which gives a low-rank decomposition of $\bm M$.
\begin{prop}\label{clm:VVT}
Let $\bm V=\left(\mathbf{1}_{1\times K}\otimes\bm R\right)\odot\left(\bm H\otimes \mathbf{1}_{d\times d}\right)\in\R^{nd\times Kd}$. Then, for the matrix $\bm M\in\R^{nd\times nd}$ defined in \eqref{eq:def-M}, $\bm M=\bm V\bm V^\top $. Moreover, $\bm V^\top\bm V=m\bm I_{Kd}$ and hence $\frac{1}{\sqrt{m}}\bm V$ is orthonormal.
\end{prop}
\begin{rmk}
Before heading for further results, let us take a closer look at the matrix $\bm V$ constructed in the Proposition. It incorporates a similar structure with $\bm H$, with exactly one nonzero block in each block row and exactly $m$ many nonzero blocks in each block column. It also encapsulates the group information $\bm R$, because $\bm R_i$ is precisely the unique nonzero block in the $i$-th block row of $\bm V$. For example, given $\bm R=\begin{pmatrix}\bm R_1\\\bm R_2\\\bm R_3\\\bm R_4\end{pmatrix}$ and $\bm H=\begin{pmatrix}1&0\\1&0\\0&1\\0&1\end{pmatrix}$, we have $\bm V=\begin{pmatrix}\bm R_1&\mathbf{0}\\\bm R_2&\mathbf{0}\\\mathbf{0}&\bm R_3\\\mathbf{0}&\bm R_4\end{pmatrix}$. Therefore, one may regard $\bm V$ as a generalization of $\bm H$ with group information $\bm R$, and $\bm H$ a degenerated case of $\bm V$ when $\G=\SO(1)=\{1\}$.
\end{rmk}
Proposition \ref{clm:VVT} then reformulates Problem \ref{pr:SDP} to a quadratic program subject to non-convex constraints, because
\begin{align}
    \max_{\substack{\bm R\in\G^n\\\bm H\in\mathcal{H}}}\left<\bm A, \bm M\right>&=\max_{\bm V\in\mathcal{E}}\left<\bm A, \bm V\bm V^\top\right>\nonumber\\
    &=\max_{\bm V\in\mathcal{E}}\mathrm{tr}(\bm V^\top \bm A\bm V),\label{pr:entangledQP-tight}\tag{\textsf{JointQP-}$\G$}
\end{align}
where the feasible region $\mathcal{E}:=\left\{\left(\mathbf{1}_{1\times K}\otimes\bm R\right)\odot\left(\bm H\otimes \mathbf{1}_{d\times d}\right)\;\middle|\;\bm R\in\G^n, \bm H\in\mathcal{H}\right\}$ is non-convex. Since $\SO(d)$ is a subgroup of $\OO(d)$, Problem \ref{pr:entangledQP-tight} can be relaxed to
\begin{equation}\label{pr:entangledQP}\tag{\textsf{JointQP-}$\OO(d)$}
    \max_{\bm V\in\mathcal{F}}\mathrm{tr}(\bm V^\top \bm A\bm V),
\end{equation}
where $\mathcal{F}:=\left\{\left(\mathbf{1}_{1\times K}\otimes\bm R\right)\odot\left(\bm H\otimes \mathbf{1}_{d\times d}\right)\setbar\bm R\in\OO(d)^n, \bm H\in\mathcal{H}\right\}$ is again a non-convex feasible region. Carefully note that the outliers in $\mathrm{SGBM}(n,p,q,K,d,\G)$ are still uniformly sampled from $\G$, although all the group elements are allowed to fall in the relaxed region $\OO(d)^n$.

When $\mathcal{E}\neq\mathcal{F}$, it is necessary to introduce a rounding function $\mathcal{R}$ that maps a matrix $\bm V\in\mathcal{F}$ back to $\mathcal{E}$. The rounding function is defined blockwise: 
\begin{equation*}
    \mathcal{R}(\bm V)_{ij}=\det(\bm V_{ij})\bm V_{ij},\;\forall i\in[n],j\in[K].
\end{equation*}
As we will prove at the end of Section \ref{subsec:proof-thm-1}, the relaxation from $\G^n$ to $\OO(d)^n$ is sufficiently tight (up to an observation-invariant transformation) thanks to this simple rounding function $\mathcal{R}$, which makes it possible to continue the study of both orthogonal and rotational scenarios in Problem \ref{pr:entangledQP}.

Finally, we say 
\begin{equation*}
\Pi_\mathcal{F}(\bm X):=\argmin_{\bm W\in\mathcal{F}}\nm{\bm W-\bm X}_F
\end{equation*}
is the projection of $\bm X$ onto $\mathcal{F}$, for arbitrary $\bm X\in\R^{nd\times Kd}$. An efficient algorithm for computing the projection will be presented in Section \ref{sec:methodology}.

\subsection{The metric of estimation error}\label{subsec:est-err}

Consider an arbitrary $\bm V\in\mathcal{E}$. Similar to the permutation invariance mentioned in Section \ref{sec:preliminaries}, any permutation of the clusters makes no difference under the settings of SGBM, since the generation of observation matrix does not depend on the specific numbering of clusters. Formally, $\bm V\bm Q\in\mathcal{E}$ would yield the same probabilistic distribution of observation as $\bm V$, where $\bm Q\in\mathcal{Q}:=\left\{\bm P\otimes\bm I_d\setbar\bm P\in\mathcal{S}_K\right\}$ and $\mathcal{S}_K$ is the permutation group on $[K]$. Moreover, right-multiplying any element of $\G$ commonly on each cluster neither affects the observation it induces. This is because for any pair of nodes $i$ and $j$ belonging to the same cluster with group elements $\bm R_i,\bm R_j$, their relative measurement always remains intact after a common right multiplication:
$$(\bm R_i\bm U)(\bm R_j\bm U)^\top=\bm R_i\bm U\bm U^\top\bm R^{\top}_j=\bm R_i\bm R^{\top}_j,\;\forall\bm U\in\G.$$
Therefore, $\bm V\bm W\in\mathcal{E}$ also yields the same probabilistic distribution of observation as $\bm V$, where
$$\bm W\in\left\{\begin{pmatrix}
    \bm U_1&&&\\
    &\bm U_2&&\\
    &&\ddots&\\
    &&&\bm U_K
    \end{pmatrix}\setbar\bm U_1,...,\bm U_K\in\G\right\}=:\mathrm{bdiag}\left(\G^K\right).
$$
One can further deduce that all the combinations of these two types of operations in $\mathcal{Q}$ and $\mathrm{bdiag}\left(\G^K\right)$ still preserves the observation, which constitute the following group:
\begin{defi}
    Let $\mathcal{P}_K(\G)$ be a group equipped with the usual matrix multiplication as the group operation, such that the matrix $\bm Q\in\R^{Kd \times Kd}$ belongs to $\mathcal{P}_K(\G)$ if and only if there exists $\bm R\in\G^K$ and $\bm P\in\mathcal{S}_K$ such that $\bm Q=\mathrm{bdiag}(\bm R)(\bm P\otimes\bm I_d)$.
\end{defi}
\begin{rmk}
    For example, given $\bm R_1,\bm R_2,\bm R_3\in\OO(d)$, the matrix $\bm Q=\begin{pmatrix}
        &&\bm R_1\\
        \bm R_2&&\\
        &\bm R_3&
    \end{pmatrix}\in\mathcal{P}_3(\OO(d))$. Note that $\mathcal{P}_K(\G)$ is a subgroup of $\OO(Kd)$. To see this, we first observe that it is a subset of $\OO(Kd)$. Then, for arbitrary $\bm Q_1,\bm Q_2\in\mathcal{P}_K(\G)$, $\bm Q_1\bm Q_2^{-1}=\bm Q_1\bm Q_2^\top\in\mathcal{P}_K(\G)$. Hence this is a subgroup of $\OO(Kd)$.
\end{rmk}

With this definition at hand, an equivalence class of $\bm V$ is immediately given by the orbit of $\bm V$ under the right group action of $\mathcal{P}_K(\G)$. To reiterate, this equivalence relation $\sim$ between two points in $\mathcal{E}$ originates from the identical probabilistic distribution of their induced observations under SGBM, which we cannot tell apart in any algorithm relying solely on these observations. It naturally induces a quotient space $\mathcal{E}/\sim$. We now define a metric for this quotient space that essentially measures the distance between the corresponding orbits of two arbitrary points.

\begin{defi}\label{def:dist}
    For any $\bm V_1,\bm V_2\in\mathcal{E}$,
    \begin{equation*}
        \mathrm{dist}_\G(\bm V_1,\bm V_2):=\min_{\bm Q\in\mathcal{P}_K(\G)}\nm{\bm V_1-\bm V_2\bm Q}_F
    \end{equation*}
    is said to be the distance between $\bm V_1$ and $\bm V_2$ in the quotient space $\mathcal{E}/\sim$.
\end{defi}
\begin{prop}
    The quotient space $\mathcal{E}/\sim$ equipped with the distance function $\mathrm{dist}$ is a metric space. 
\end{prop}

\begin{proof}
    To establish triangle inequality: for arbitrary $\bm V_1,\bm V_2,\bm V_3\in\mathcal{E}$, let 
    $$
    \bm Q_{12}=\argmin_{\bm Q\in\mathcal{P}_K(\G)}\nm{\bm V_1-\bm V_2\bm Q}_F,\;\bm Q_{23}=\argmin_{\bm Q\in\mathcal{P}_K(\G)}\nm{\bm V_2-\bm V_3\bm Q}_F.
    $$
    Then, since $\mathcal{P}_K(\G)$ is a subgroup of $\OO(Kd)$,
    \begin{align}
    &\mathrm{dist}_\G(\bm V_1,\bm V_2)+\mathrm{dist}_\G(\bm V_2,\bm V_3)=\nm{\bm V_1-\bm V_2\bm Q_{12}}_F+\nm{\bm V_2-\bm V_3\bm Q_{23}}_F\nonumber\\
    =&\nm{\bm V_1-\bm V_2\bm Q_{12}}_F+\nm{\bm V_2\bm Q_{12}-\bm V_3\bm Q_{23}\bm Q_{12}}_F\geq\nm{\bm V_1-\bm V_3\bm Q_{23}\bm Q_{12}}_F\nonumber\\
    \geq &\min_{\bm Q\in\mathcal{P}_K(\G)}\nm{\bm V_1-\bm V_3\bm Q}_F=\mathrm{dist}_\G(\bm V_1,\bm V_3).\nonumber
    \end{align}
    The other two axioms for the metric space are easy to verify.
\end{proof}

As a major motivation of Definition \ref{def:dist}, we are specially interested in properly measuring the estimation error of an arbitrary $\bm V\in\mathcal{E}$ and the ground truth $\bm V^*$. Via the distance in the quotient space, Definition \ref{def:dist} immediately gives a reasonable metric of estimation error of $\bm V$:

\begin{defi}\label{def:errorV}
Let $\bm V,\bm V^*\in\mathcal{E}$ be given and $\bm V^*$ is the ground truth. Then, $$\epsilon_\G(\bm V):=\mathrm{dist}_\G(\bm V,\bm V^*)=\min_{\bm Q\in\mathcal{P}_K(\G)}\nm{\bm V-\bm V^*\bm Q}_F$$ is said to be the estimation error of $\bm V$.
\end{defi}

\section{A non-convex methodology}\label{sec:methodology}
We now propose a \textit{Generalized Power Method} (GPM) to tackle Problem \ref{pr:entangledQP}. The approach to initializing $\bm V^0$ is given by Algorithm \ref{alg:spec-init} to be covered in detail in Section \ref{sec:spec-init}.
\begin{algorithm}[H]
    \caption{GPM} \label{alg:entGPM}
    \begin{algorithmic}[1]
    \State \textbf{Input:} the observation matrix $\bm A$, the initial point $\bm V^{0}$.
    \For {$t = 0,1,2,\dots, T-1$}
    \State $\bm V^{t+1}\leftarrow\Pi_\F(\bm A\bm V^t)$
    \EndFor
    \If {$\G=\SO(d)$}
    \State {$\bm V^T\leftarrow\mathcal{R}(\bm V^T)$}
    \EndIf
    \State \textbf{Return:} $\bm V^T$
    \end{algorithmic}
\end{algorithm}

It remains to design an algorithm that efficiently solves the projection onto $\mathcal{F}$. To this end, we first define a mapping $\mu:\R^{nd\times Kd}\to\R^{n\times K}$ that computes the blockwise nuclear norm
\begin{equation*}
    \mu(\bm X)_{ij}=\nm{\bm X_{ij}}_*=\sum_{k=1}^d\sigma_k(\bm X_{ij}),\;\forall i\in[n],j\in[K],
\end{equation*}
and present the algorithm as follows:
\begin{algorithm}[H]
    \caption{Projection onto $\mathcal{F}$} \label{alg:proj}
    \begin{algorithmic}[1]
    \State \textbf{Input:} $\bm X\in\R^{nd\times Kd}$
    \State $\bm H\leftarrow\Pi_\mathcal{H}(\mu(\bm X))$
    \State generate a sequence $\{e_i\}_{i=1}^n$ such that $\bm H_{ie_i}=1$
    \For{$i=1,2,\cdots,n$}
        \State $\bm R_{i}\leftarrow\Pi_{\OO(d)}(\bm X_{ie_i})$
    \EndFor
    \State $\bm V\leftarrow \left(\mathbf{1}_{1\times K}\otimes \bm R\right)\odot\left(\bm H\otimes \mathbf{1}_{d\times d}\right)$
    \State \textbf{Return:} $\bm V$
    \end{algorithmic}
\end{algorithm}
\begin{prop}\label{prop:proj_F-runtime}
Given that $K,d=\Theta(1)$, Algorithm \ref{alg:proj} exactly solves $\Pi_\F(\bm X)$ in $O(n\log n)$ time.
\end{prop}

As its name suggests, the design of GPM is inspired by the classical \textit{power method} as well as its variant \textit{orthogonal iteration method} \cite{golub2012matrix} that computes the dominant eigenvector and the invariant (orthogonal) subspaces of a matrix, respectively. In fact, one may have already observed resemblances between Problem \ref{pr:entangledQP} and the classical eigenvalue problems above. For example, both maximize quadratic objectives, and both the optimization variables are subject to norm and orthogonality constraints by noting that $\frac{1}{\sqrt{m}}\bm V$ is orthonormal. The algorithm in general is concise: it refines the initial guess iteratively by simply taking matrix multiplication and projection onto the relaxed feasible region $\mathcal{F}$, involving no hyperparameters. Therefore, to implement our GPM algorithm is practically convenient.

\section{Main result}\label{sec:main-result}
\newcommand{\tauBernBound}{\sqrt{2K\beta}}
\newcommand{\tauBernBoundSq}{2K\beta}

In this section, we first formally state the main result, Theorem \ref{thm:main}--\ref{thm:cond-iii}, on the guarantee of linear convergence of GPM under the metric of the estimation error per iteration in the logarithmic sparsity region of SGBM, \ie, $p,q=O\left(\frac{\log n}{n}\right)$, and then present the proof of Theorem \ref{thm:main}. We defer the proof of Theorem \ref{thm:cond-i-ii} to the Appendix, and Theorem \ref{thm:cond-iii} to section \ref{sec:spec-init} where we provide a randomized spectral clustering method for initialization of GPM. 

\subsection{Statement of the main result}
\begin{defi}
    Suppose $\bm V\in\mathcal{F}$ with a clustering function $C$, and the observation matrix $\bm A\in\R^{nd\times nd}$ is given. Let $\bm M=\mu(\bm A\bm V)$. Then the observation matrix $\bm A$ is said to preserve $\bm V$ by $\delta$-separation if $\bm M_{iC(i)}-\bm M_{ij}\geq\delta>0$ for all $i\in[n]$ and $j\neq C(i), j\in[K]$.
\end{defi}

\begin{thm}[Main]\label{thm:main}
    Suppose that the observation matrix $\bm A$ is generated by 
    $$\mathrm{SGBM}\left(n,\frac{\alpha\log n}{n},\frac{\beta\log n}{n},K,d,\G\right),$$
    where $\alpha,\beta,K,d=\Theta(1)$, and $\G\in\{\OO(d),\SO(d)\}$. Let $\bm V^0\in\mathcal{F}$ and $\bm A$ be the input of Algorithm \ref{alg:entGPM}. Let $\bm V^*\in\mathcal{E}\subset\mathcal{F}$ be the ground truth matrix that incorporates the clustering function $C^*$. Then, Algorithm \ref{alg:entGPM} outputs $\bm V^T\in\mathcal{E}$ such that $\epsilon_\G(\bm V^T)\leq\tau$ within $O\left(n\log n\log\frac{n}{\tau}\right)$ time, if the following hold:\\
    (\textit{i}) there exists a constant $\chi>0$ such that $\bm A$ perserves $\bm V^*$ by $\chi mp$-separation;\\
    (\textit{ii}) for all $i\in[n]$, $\mu(\bm A\bm V^*)_{iC^*(i)}\geq\frac{\tauBernBound}{\alpha}mp$;\\
    (\textit{iii}) for $\chi>0$ satisfying (\textit{i}), there exists 
    $$\rho> 2\sqrt{2}\sqrt{\frac{d^2}{\chi^2}+\frac{\alpha^2}{2K\beta}}$$
    such that $\epsilon_{\OO(d)}(\bm V^0)\leq\frac{\sqrt{m}}{\rho}$.
\end{thm}
\begin{rmk}
    Note that Theorem \ref{thm:main} is a \textit{deterministic} statement given the three conditions above. From a qualitative perspective, the first and second conditions describe a desirable interaction between the observation $\bm A$ and the ground truth $\bm V^*$ that sufficiently separates the mutual measurement from noise. As we would see through the proof of the theorem, they provide a decent bound for the Lipschitz constant in our analysis of $\Pi_\mathcal{F}$. The third condition, on the other hand, demands a reasonable initial guess $\bm V^0$ of $O(\sqrt{m})$ error before the refinement steps in GPM. Since the largest possible estimation error $$\sup_{\bm V\in\mathcal{F}}\epsilon_{\OO(d)}(\bm V)=\sqrt{2Kdm}=O(\sqrt{m})$$
    is of the same order, condition (\textit{iii}) is rather tolerant.
\end{rmk}

For a sufficiently large $n$, the following two companion theorems provide a probablistic lower bound for the parameters $(\alpha,\beta,K,d)$ such that condition (\textit{i}) and (\textit{ii}) hold with probability at least $1-n^{-\Omega(1)}$, and condition (\textit{iii}) holds with probability at least $1-\left(\log n\right)^{-\Omega(1)}$.

\begin{thm}[Condition (\textit{i}) and (\textit{ii})]\label{thm:cond-i-ii}
    Suppose that $\alpha,\beta,K,d=\Theta(1)$, and the observation matrix $\bm A$ is generated by $\mathrm{SGBM}$ for the given ground truth matrix $\bm V^*$ with the clustering mapping $C^*$. If
    \begin{numcases}{}
        &$\tauBernBound<\alpha$,\label{eq:region-for-i-ii-cond1}\\
        &$\alpha-\tauBernBound\log\frac{e\alpha}{\tauBernBound}>K$,\label{eq:region-for-i-ii-cond2}
     \end{numcases}
    \noindent then condition (\textit{i}) and (\textit{ii}) in Theorem \ref{thm:main} happen simultaneously with probability at least $1-n^{-\Omega(1)}$ for a sufficiently large $n$.
\end{thm}

\begin{thm}[Condition (\textit{iii})]\label{thm:cond-iii}
    Suppose that $\alpha,\beta,K,d=\Theta(1)$, and the observation matrix $\bm A$ is generated by $\mathrm{SGBM}$ for the given ground truth matrix $\bm V^*$ with the clustering mapping $C^*$. Then, there exists an algorithm that generates an initial $\bm V^0$ satisfying condition (\textit{iii}) in Theorem \ref{thm:main} with proability at least $1-(\log n)^{-\Omega(1)}$ for a sufficiently large $n$.
\end{thm}

\subsection{Proof of Theorem \ref{thm:main}}\label{subsec:proof-thm-1}
To prove Theorem \ref{thm:main}, we first show the linear convergence of Algorithm \ref{alg:entGPM} in the quotient space $\mathcal{F}/\sim$ under the relaxed metric of estimation error, $\epsilon_{\OO(d)}$. Then, we present the tightness of the rounding procedure $\mathcal{R}$ that brings forth Theorem \ref{thm:main} as a direct consequence.

\subsubsection*{Linear convergence under the metric $\epsilon_{\OO(d)}$}
This section is devoted to proving the following theorem on the linear convergence of Algorithm \ref{alg:entGPM} under the relaxed metric of estimation error $\epsilon_{\OO(d)}$.

\begin{thm}\label{thm:linear-conv}
    Suppose that the observation matrix $\bm A$ is generated by 
    $$\mathrm{SGBM}\left(n,\frac{\alpha\log n}{n},\frac{\beta\log n}{n},K,d,\G\right),$$
    where $\alpha,\beta,K,d=\Theta(1)$, and $\G\in\{\OO(d),\SO(d)\}$. Let $\bm V^0\in\mathcal{F}$ and $\bm A$ be the input of Algorithm \ref{alg:entGPM}, and $\{\bm V^1,\bm V^2,...\}$ a sequence generated by the iterations in Algorithm \ref{alg:entGPM}. Let $\bm V^*\in\mathcal{E}\subset\mathcal{F}$ be the ground truth matrix that incorporates the clustering function $C^*$. Then,
    $$\epsilon_{\OO(d)}(\bm V^{t+1})\leq\frac{1}{2}\epsilon_{\OO(d)}(\bm V^t)$$
    for arbitrary non-negative integer $t$, if condition (\textit{i})--(\textit{iii}) in Theorem \ref{thm:main} hold.
\end{thm}

Recall that GPM iteratively updates the variable by the map $\Pi_\mathcal{F}(\bm A\;\cdot\;):\mathcal{F}\to\mathcal{F}$. Therefore, in order to prove Theorem \ref{thm:linear-conv}, it should be crucial to identify the important properties of the projection operator $\Pi_\mathcal{F}$. We first make two useful observations in the following Lemmas:

\begin{lem}\label{lem:Pi-XQ}
    For any $\bm Q\in\mathcal{P}_K(\OO(d))$ and $\bm X\in\R^{nd\times Kd}$, $\Pi_{\mathcal{F}}(\bm X\bm Q)=\Pi_\mathcal{F}(\bm X)\bm Q$.
\end{lem}
\begin{lem}\label{lem:fixed-point}
    If condition (\text{i}) in Theorem \ref{thm:main} holds, then $\Pi_\F(\bm A\bm V^*)=\bm V^*$.
\end{lem}

The core of the proof lies in controlling the behavior of the projection operator $\Pi_\mathcal{F}$ so that the estimation error after each update can be bounded. One possible way, as Proposition \ref{thm:Lipschitz} follows, is to show $\Pi_\mathcal{F}$ possesses a Lipschitz-like property of linearly controlling the Frobenius distance of two points after a projection.
\begin{prop}\label{thm:Lipschitz}
    Let $\bm X=\bm A\bm V^*$. If condition (\text{i}) and (\text{ii}) in Theorem \ref{thm:main} hold, then
    $$
    \nm{\Pi_\mathcal{F}(\bm X)-\Pi_\mathcal{F}(\bm X')}_F\leq\frac{2}{mp}\sqrt{\frac{d^2}{\chi^2}+\frac{\alpha^2}{\tauBernBoundSq}}\nm{\bm X-\bm X'}_F.
    $$
    for any $\bm X'\in\R^{nd\times Kd}$.
\end{prop}

For the sequence $\{\bm V^0,\bm V^1,\bm V^2,...\}$ generated by GPM, denote $\bm Q^t=\argmin_{\bm Q\in\mathcal{P}_K(\G)}\nm{\bm V-\bm V^*\bm Q}_F$ for all $t\geq 0$. Given that condtion (\textit{i}) and (\textit{ii}) in Theorem \ref{thm:main} hold, we have
\begin{align}
    \nm{\bm V^{t+1}-\bm V^*\bm Q^{t+1}}_F&\leq\nm{\Pi_\F(\bm A\bm V^t)-\bm V^*\bm Q^t}_F=\nm{\Pi_\F(\bm A\bm V^t\bm Q^{t\top})-\bm V^*}_F\nonumber\\
    &=\nm{\Pi_\F(\bm A\bm V^t\bm Q^{t\top})-\Pi_\F(\bm A\bm V^*)}_F\leq\frac{2}{mp}\sqrt{\frac{d^2}{\chi^2}+\frac{\alpha^2}{2K\beta}}\nm{\bm A\bm V^t\bm Q^{t\top}-\bm A\bm V^*}_F\nonumber\\
    &=\frac{2}{mp}\sqrt{\frac{d^2}{\chi^2}+\frac{\alpha^2}{2K\beta}}\nm{\bm A(\bm V^t-\bm V^*\bm Q^t)}_F,\nonumber
\end{align}
where Lemma \ref{lem:Pi-XQ} yields the first equality, Lemma \ref{lem:fixed-point} yields the second equality, and the second inequality is due to Proposition \ref{thm:Lipschitz}. We can proceed to obtain 
\begin{align}
&\nm{\bm A(\bm V^t-\bm V^*\bm Q^t)}_F\leq\nm{(\bm A-p\bm V^*\bm V^{*\top})(\bm V^t-\bm V^*\bm Q^t)}_F+\nm{p\bm V^*\bm V^{*\top}(\bm V^t-\bm V^*\bm Q^t)}_F\nonumber\\
\leq&\nm{\bm A-p\bm V^*\bm V^{*\top}}\nm{\bm V^t-\bm V^*\bm Q^t}_F+\nm{p\bm V^*\bm V^{*\top}(\bm V^t-\bm V^*\bm Q^t)}_F\nonumber\\
=&\nm{\bm A-p\bm V^*\bm V^{*\top}}\nm{\bm V^t-\bm V^*\bm Q^t}_F+m\sqrt{m}p\nm{\frac{1}{m}(\bm V^*\bm Q^t)^\top \bm V^t-\bm I_{Kd}}_F.\label{eq:intermediate}
\end{align}
Since we are expecting some relationship between $\nm{\bm V^{t+1}-\bm V^*\bm Q^{t+1}}_F$ and $\nm{\bm V^t-\bm V^*\bm Q^t}_F$, it remains to bound the two components $\nm{\bm A-p\bm V^*\bm V^{*\top}}$ and $\nm{\frac{1}{m}\left(\bm V^*\bm Q^t\right)^\top\bm V^t-\bm I_{Kd}}_F$ respectively, as is presented in the following two Propositions.

\begin{prop}\label{prop:A-pVVT}
    There exists $c_1, c_2, c_3>0$ such that
    \begin{equation*}
        \nm{\bm{A}-p\bm{V}^*\bm{V}^{*\top}}\leq c_1\sqrt{qm}+c_2\sqrt{pm}+c_3\sqrt{\log n}
    \end{equation*}
    with probability at least $1-n^{-\Omega(1)}$.
\end{prop}
\begin{rmk}
    In the logarithmic sparsity region, Proposition \ref{prop:A-pVVT} gives a $O(\sqrt{\log n})$ bound. Different from the techniques in Theorem 6 of \cite{liu2020unified} for controlling a similar term, the result here does not rely on the celebrated matrix Bernstein inequality \cite{Tropp2011matrixbound}, which loosely controls the term at $O(\log n)$.
\end{rmk}

\begin{prop}\label{prop:Z-I-bound}
    For $\rho> 0$, 
    $$m\nm{\frac{1}{m}\left(\bm V^*\bm Q^t\right)^{\top}\bm V^t-\bm I_{Kd}}_F\leq \frac{1}{2}\sqrt{1+\frac{1}{d}}\frac{\sqrt{m}}{\rho}\nm{\bm V^t-\bm V^*\bm Q^t}_F$$
    if $\nm{\bm V^t-\bm V^*\bm Q^t}_F\leq\frac{\sqrt{m}}{\rho}$.
\end{prop}

\begin{proof}[Proof for Theorem \ref{thm:linear-conv}]
    Invoking Proposition \ref{prop:A-pVVT} and \ref{prop:Z-I-bound}, \eqref{eq:intermediate} becomes
    \begin{align}
        \nm{\bm V^{t+1}-\bm V^*\bm Q^{t+1}}_F&\leq\frac{2}{mp}\sqrt{\frac{d^2}{\chi^2}+\frac{\alpha^2}{2K\beta}}\left(\frac{1}{2}\sqrt{1+\frac{1}{d}}\frac{mp}{\rho}+c_0\sqrt{\log n}\right)\nm{\bm V^t-\bm V^*\bm Q^t}_F\nonumber\\
        &\leq\frac{\sqrt{2}+o(1)}{\rho}\sqrt{\frac{d^2}{\chi^2}+\frac{\alpha^2}{2K\beta}}\nm{\bm V^t-\bm V^*\bm Q^t}_F\nonumber
    \end{align}
    for a sufficiently large $n$. When 
    $$\rho> 2\sqrt{2}\sqrt{\frac{d^2}{\chi^2}+\frac{\alpha^2}{2K\beta}},$$
    we have $\epsilon_{\OO(d)}(\bm V^{t+1})\leq\frac{1}{2}\epsilon_{\OO(d)}(\bm V^t)$. When condition (\textit{iii}) in Theorem \ref{thm:main} holds, the linear decay of $\epsilon(\bm V^t)$ inductively applies to arbitrary nonnegative integer $t$, which establishes Theorem \ref{thm:linear-conv}.
\end{proof}

\subsubsection*{Tightness of rounding}
It is an immediate implication of Theorem \ref{thm:linear-conv} that, prior to the rounding procedure, Algorithm \ref{alg:entGPM} obtains $\bm V^T$ such that $\epsilon_{\OO(d)}(\bm V^T)\leq\tau$ within $T=O\left(\log\frac{n}{\tau}\right)$ iterations. In each iteration of GPM, the projection takes $O(n\log n)$ time according to Proposition \ref{prop:proj_F-runtime}, and the time complexity of matrix multiplication can also achieve $O(n\log n)$ due to their sparse structures. Therefore, the total time complexity to obtain such a $\bm V^T$ is $O\left(n\log n\log \frac{n}{\tau}\right)$. Compared with Theorem \ref{thm:main}, the very difference lies in the metric of estimation error $\epsilon_\G$ after the rounding procedure $\mathcal{R}$. We handle the tightness of $\mathcal{R}$, hence the tightness of relaxation from $\mathcal{E}$ to $\mathcal{F}$ unresolved in Section \ref{sec:formulation}, by the following Proposition.
\begin{prop}\label{prop:tightness}
    Suppose that $\G=\SO(d)$, and $\bm V^*\in\mathcal{E}$ is the ground truth. Then, for any $\bm V\in\mathcal{F}$ such that $\epsilon_{\OO(d)}(\bm V)<\sqrt{2}$, $\epsilon_{\G}(\mathcal{R}(\bm {V}))=\epsilon_{\OO(d)}(\bm V)$.
\end{prop}

In theory and practice, the tolerance constant $\tau$ is usually far below $\sqrt{2}$. Therefore, combined with Theorem \ref{thm:linear-conv}, Proposition \ref{prop:tightness} completes the proof of Theorem \ref{thm:main}.

\section{Randomized spectral clustering}\label{sec:spec-init}
We now propose a randomized spectral clustering algorithm that generates an initial point $\bm V^0$ satisfying condition (\textit{iii}) in Theorem \ref{thm:main} with high probability. Therefore, this algorithm is able to serve as a qualified initializer of GPM. 

\begin{algorithm}[H]
    \caption{Randomized spectral clustering} \label{alg:spec-init}
    \begin{algorithmic}[1]
    \State \textbf{Input:} the observation matrix $\bm A$
    \State \textbf{Initialize: } $\bm R^0\in\R^{nd\times d}$
    \State generate $\bm H^0$ by Algorithm 2 in \cite{gao2015achieving}
    \State $\bm H^0\leftarrow\Pi_\mathcal{H}(\bm H^0)$
    \State $\hU\leftarrow\mathrm{eigs}_{Kd}(\bm A)$
    \For {$i\in[K]$}
    \State pick $\tau_i\in\mathcal{I}^0_i$ uniformly randomly \texttt{//randomized pivot selection}
    \For {$v\in\mathcal{I}_i^0$}
        \State $\bm R^0_v\leftarrow\argmin_{\bm R\in\mathcal{O}(d)}\nm{\hU_{v\times}-\bm R\hU_{\tau_i\times}}_F=\Pi_{\mathcal{O}(d)}\left(\hU_{v\times}\hU_{\tau_i\times}^\top \right)$
    \EndFor
    \EndFor
    \State $\bm V^0\leftarrow \left(\mathbf{1}_{1\times K}\otimes \bm R^{0}\right)\odot\left(\bm H^0\otimes\mathbf{1}_{d\times d}\right)$
    \State \textbf{Return:} $\bm V^0$
    \end{algorithmic}
\end{algorithm}

While the community structure and group information are jointly optimized in GPM for exact recovery, this initialization algorithm adopts an intuitive two-stage design as condition (\textit{iii}) in Theorem \ref{thm:main} tolerates a rough estimation. The first stage generates a preliminary guess of the community structure $\bm H$. Similar to \cite{wang2021optimal}, we invoke Algorithm 2 in \cite{gao2015achieving}, a greedy spectral clustering method, to obtain a (imbalanced) clustering matrix, which is then rounded to a balanced clustering matrix $\bm H^0\in\mathcal{H}$. Theory developed in \cite{wang2021optimal} has shown that for any numerical constant $C$,
\begin{equation}\label{eq:H0-err}
    \epsilon(\bm H^0)\leq \sqrt{\frac{Cn}{\log n}}
\end{equation}
with probability at least $1-n^{-\Omega(1)}$. After obtaining an initialization $\bm H$, the group elements are estimated in the second stage. For each hypothesized cluster $i$ obtained previously, a pivot node denoted by $\tau_i$ is located by a randomized mechanism. Then, the relative group transformation between the pivot and any other node is estimated utilizing the corresponding block rows of $\hU=\mathrm{eigs}_{Kd}(\bm A)$. We consolidate the two estimations into a matrix $\bm V^0\in\mathcal{F}$ as the return of the algorithm.

As the following Proposition points out, the algorithm has a decent theoretical performance regarding the estimation precision of $\bm V^0$ it generates. 
\begin{prop}\label{prop:rand-spec-init}
    Suppose that $\alpha,\beta,K,d=\Theta(1)$, and the observation matrix $\bm A$ is generated by $\mathrm{SGBM}$ for the given ground truth matrix $\bm V^*$. Then, for any given constant $\rho>0$, Algorithm \ref{alg:spec-init} generates a $\bm V^0$ such that
    $$\epsilon_{\OO(d)}\left(\bm V^0\right)\leq \frac{\sqrt m}{\rho}$$
    with probability at least $1-(\log n)^{-\Omega(1)}$ for a sufficiently large $n$.
\end{prop}
\begin{proof}[Proof of Theorem \ref{thm:cond-iii}]
In fact, equipped with Proposition \ref{prop:rand-spec-init}, Theorem \ref{thm:cond-iii} is immediately established. However, we remark that this is not necessarily the unique algorithm fulfilling the requirements therein.
\end{proof}

\section{Numerical experiments}\label{sec:numerics}
In order to corroborate the theoretical analysis completed previously, this section studies and evaluates the performance of GPM through different numerical experiments including its phase transition behavior, convergence performance, and CPU time. When necessary, a comparison is also drawn between GPM and the SDP method \cite{fan2021joint}. All the simulations are conducted via MATLAB R2021a on a workstation hosting a 64-bit Windows 10 environment with 256GB RAM and Intel(R) Xeon(R) CPU E5-2699 2.20GHz 2-processor CPU.

\subsection{Phase transition}
\begin{figure}[h!]
    \centering
    \begin{subfigure}[b]{0.325\textwidth}
        \includegraphics[trim = 45mm 80mm 45mm 90mm, clip, width=\textwidth]{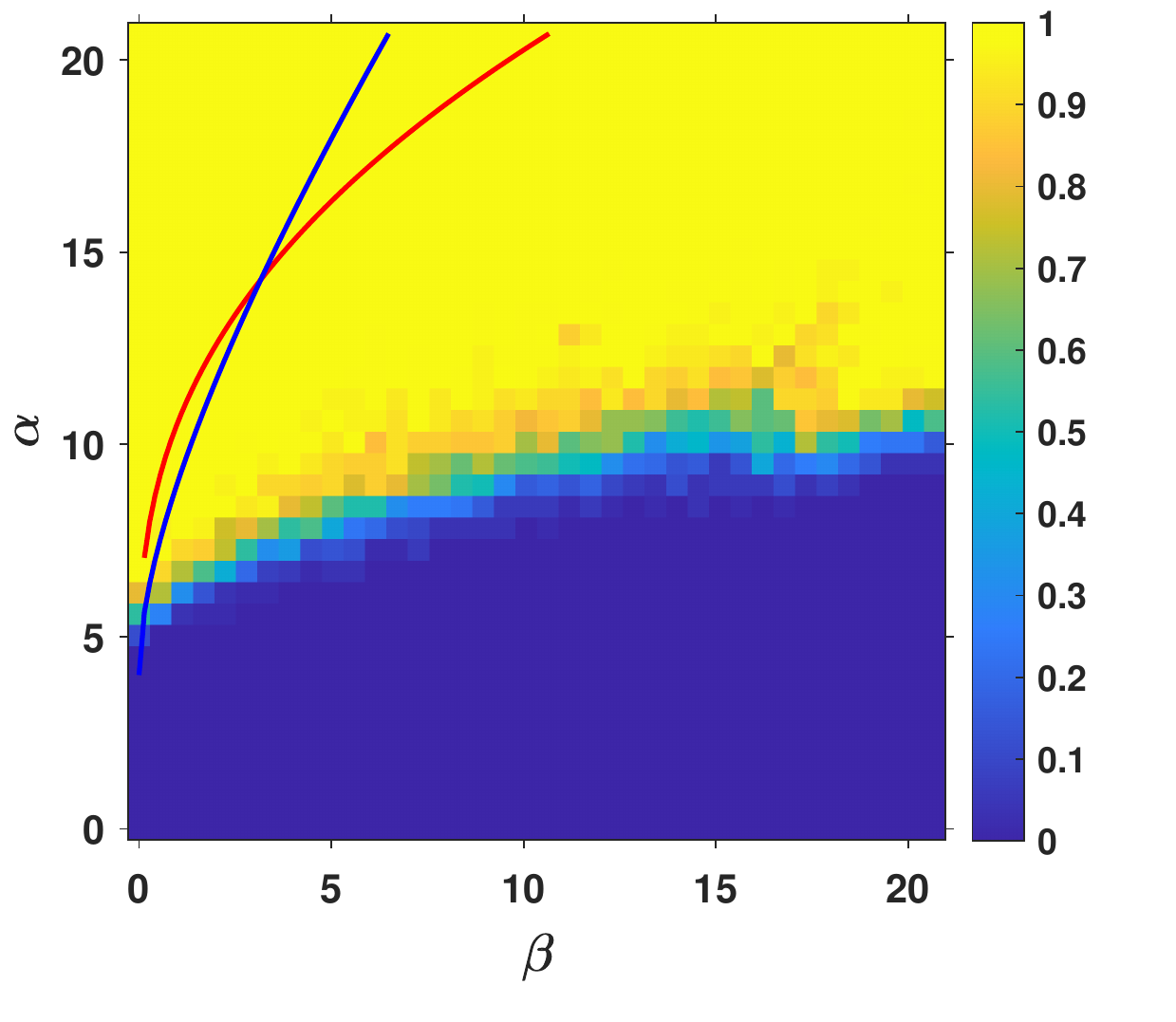}\label{fig:ps_a}
        \caption{$\G=\OO(3), n=100, K=4$}
    \end{subfigure}
    \hfill
    \begin{subfigure}[b]{0.325\textwidth}
        \includegraphics[trim = 45mm 80mm 45mm 90mm, clip, width=\textwidth]{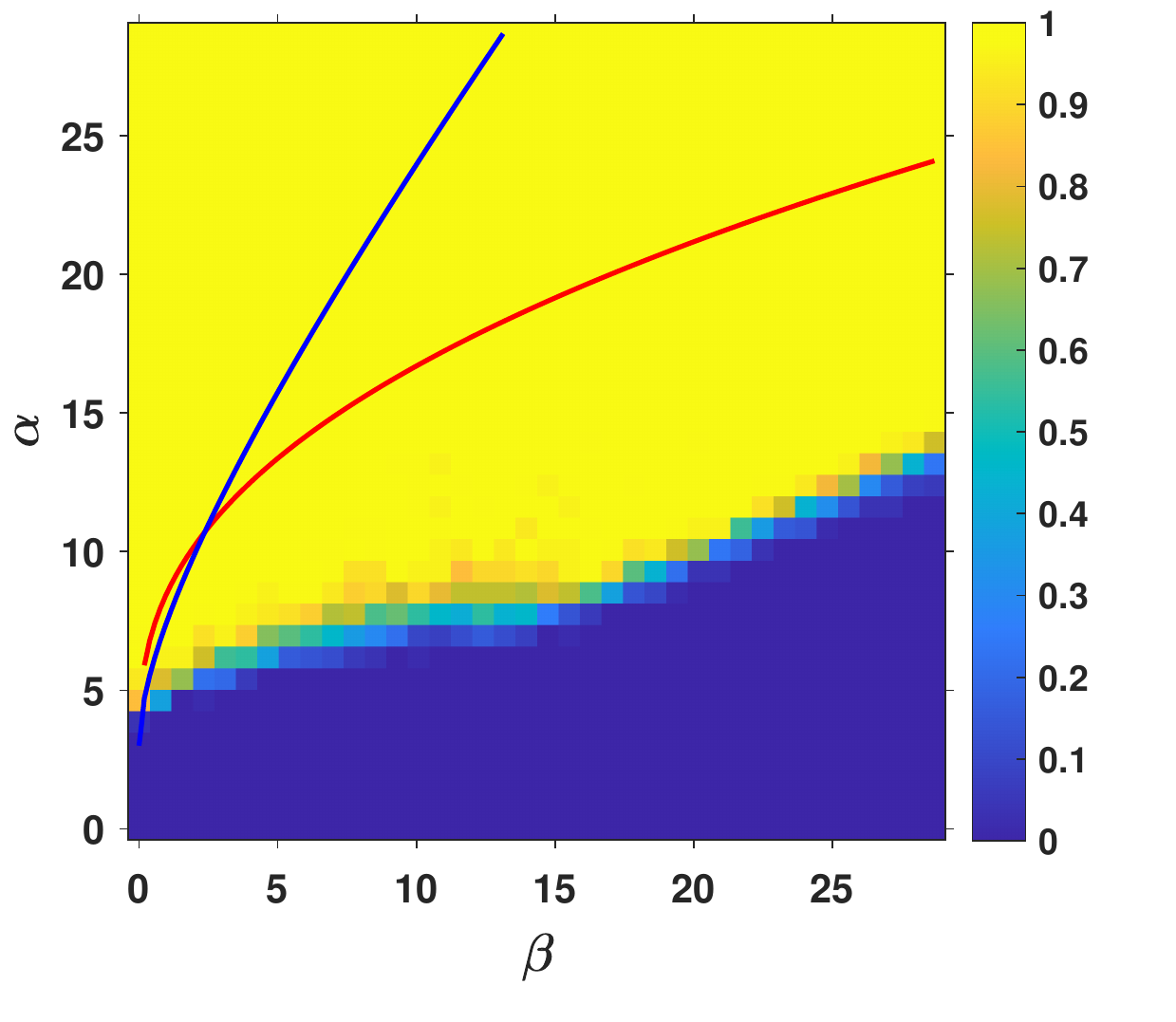}\label{fig:ps_b}
        \caption{$\G=\OO(3), n=150, K=3$}
    \end{subfigure}
    \hfill
    \begin{subfigure}[b]{0.325\textwidth}
        \includegraphics[trim = 45mm 80mm 45mm 90mm, clip, width=\textwidth]{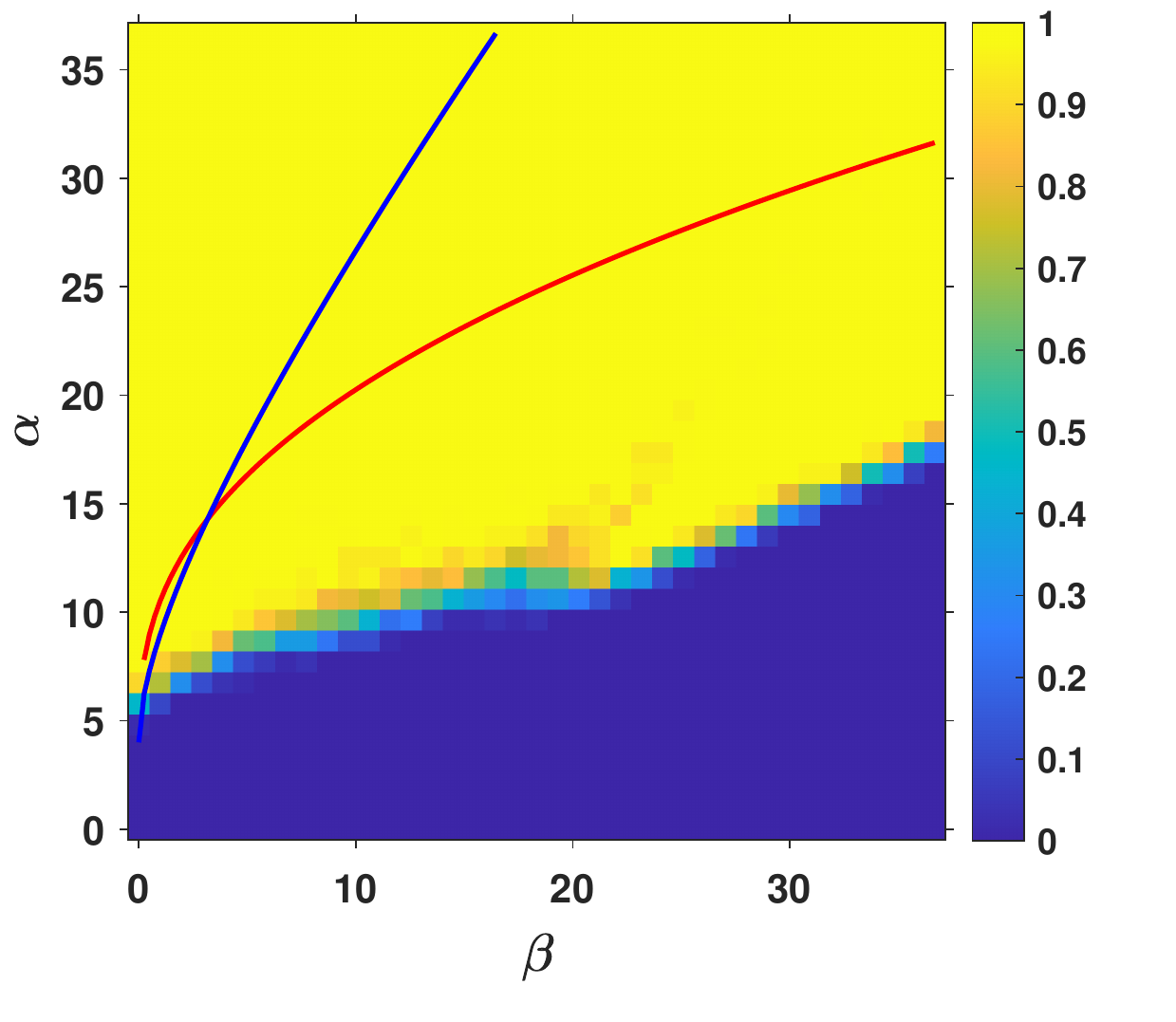}\label{fig:ps_c}
        \caption{$\G=\OO(3), n=200, K=4$}
    \end{subfigure}\\
    
    \begin{subfigure}[b]{0.325\textwidth}
        \includegraphics[trim = 45mm 80mm 45mm 90mm, clip, width=\textwidth]{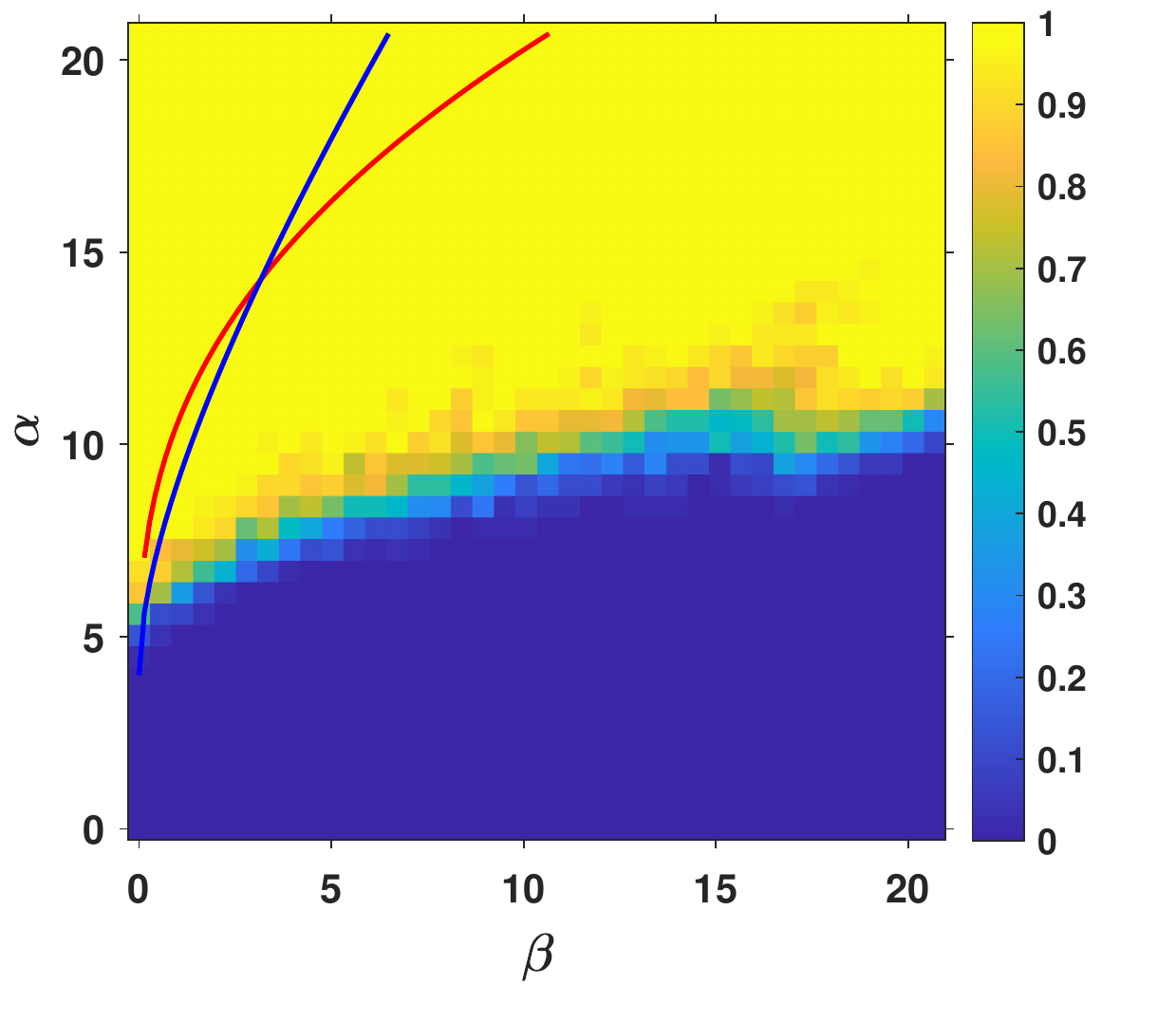}\label{fig:ps_d}
        \caption{$\G=\SO(3), n=100, K=4$}
    \end{subfigure}
    \hfill
    \begin{subfigure}[b]{0.325\textwidth}
        \includegraphics[trim = 45mm 80mm 45mm 90mm, clip, width=\textwidth]{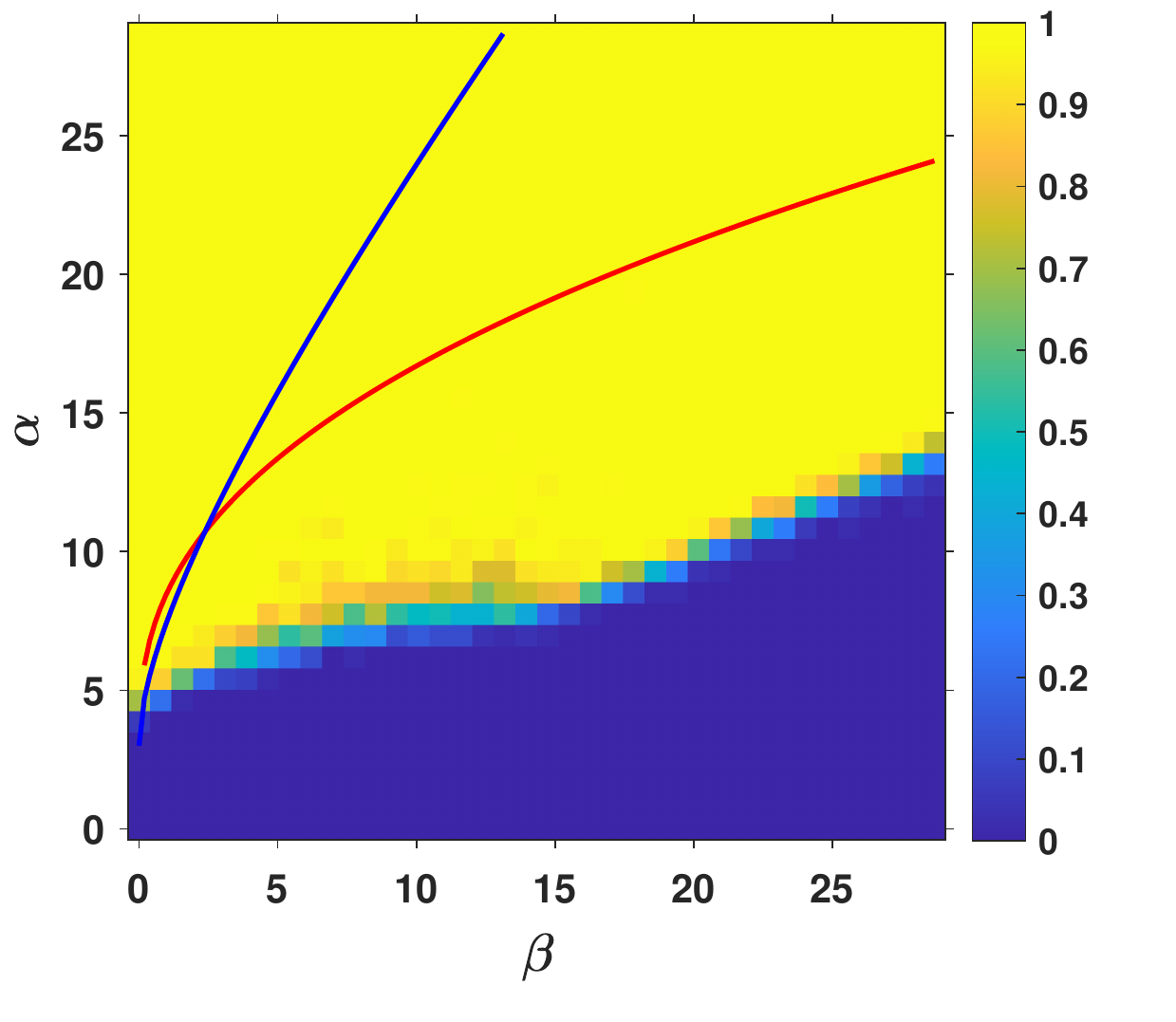}\label{fig:ps_e}
        \caption{$\G=\SO(3), n=150, K=3$}
    \end{subfigure}
    \hfill
    \begin{subfigure}[b]{0.325\textwidth}
        \includegraphics[trim = 45mm 80mm 45mm 90mm, clip, width=\textwidth]{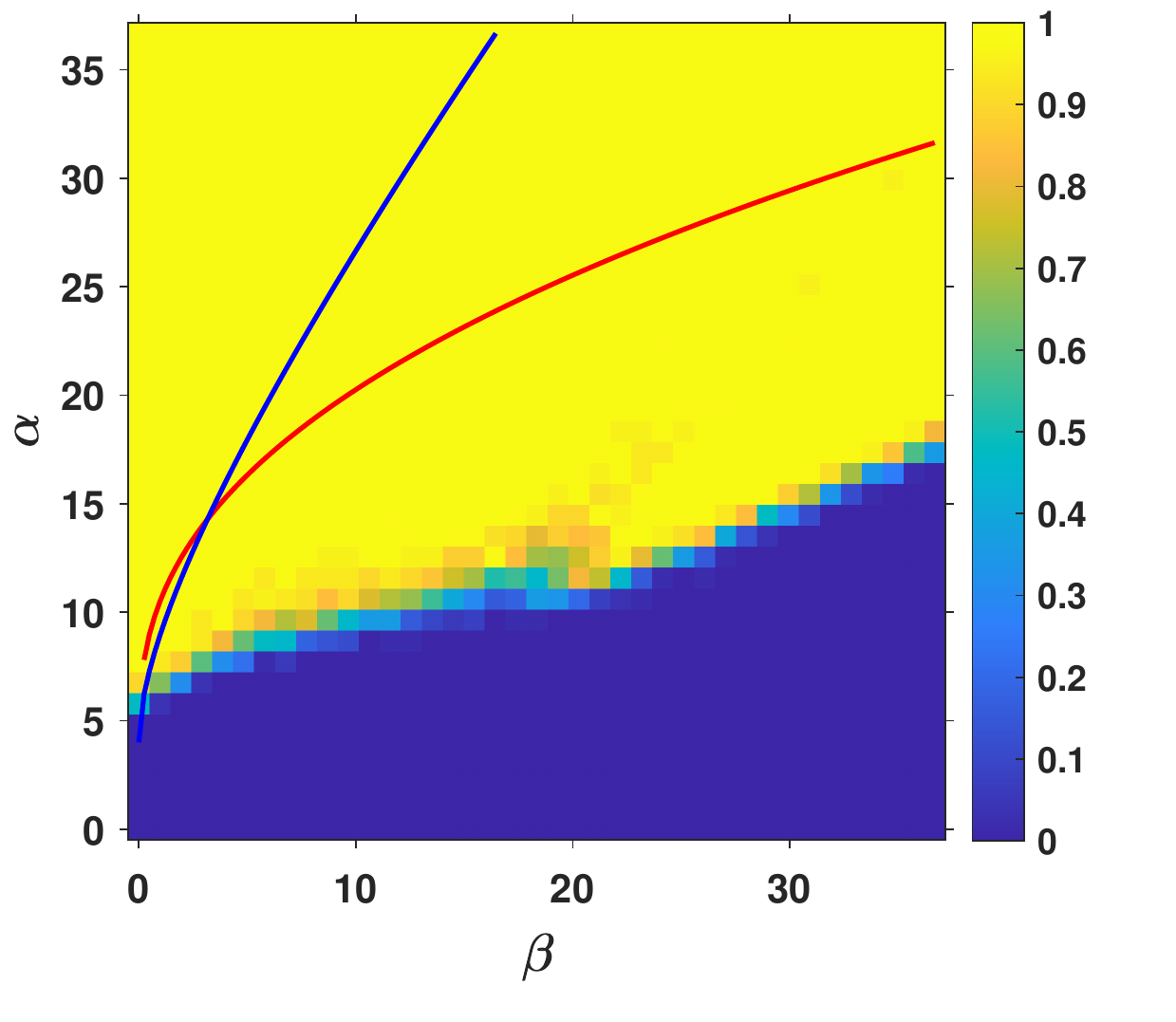}\label{fig:ps_f}
        \caption{$\G=\SO(3), n=200, K=4$}
    \end{subfigure}
    \caption{Phase transition results on GPM with three different pairs of parameters $(n,K)=(100,4)$, $(150,3)$, and $(200,4)$ in both orthogonal and rotational scenarios. The theoretical threshold for pure community detection $\sqrt{\alpha}-\sqrt{\beta}=\sqrt{K}$ is plotted in blue; the improved lower bound claimed in Theorem \ref{thm:cond-i-ii} is plotted in red.}
    \label{fig:phase-transition}
\end{figure}

We first report the phase transition behavior of GPM for both orthogonal and rotational scenarios. For each selected pair of parameters $(n,K)$, we increment $\alpha$ and $\beta$ from $0$ to $\frac{n}{\log n}$, and generate the observation under $\mathrm{SGBM}(n,\alpha\log n/n,\beta\log n/n, K, d, \G)$ mechanism within each iteration. Then GPM is invoked $N=50$ times to attempt to recover the ground truth. We regard an attempt successful if and only if $\epsilon_\G(\bm V^\top)\leq\tau=10^{-3}$, and the \textit{rate of success} given $\alpha$ and $\beta$ is thereby defined to be
\begin{equation*}
    r(\alpha,\beta)=\frac{\text{number of successes}}{N}.
\end{equation*}
We plot $r(\alpha,\beta)$ versus the change of $\alpha$ and $\beta$ in Figure \ref{fig:phase-transition}, together with the theoretical threshold for pure community detection $\sqrt{\alpha}-\sqrt{\beta}=\sqrt{K}$ (blue) and the lower bound claimed in Theorem \ref{thm:cond-i-ii} (red). In all the experiments, the results clearly exhibit a behavior of phase transition, and the region of failure is sufficiently controlled by the lower bound. The gap inbetween indicates that even the improved lower bound is still not tight for GPM, which calls for further analysis of the algorithm. One can also observe that the phase transition behave slightly differently for small and large settings of $\alpha$ and $\beta$: empirically, the boundary delimiting failure and success concaves down for smaller $\alpha$ and $\beta$, while it grows near-linear for larger parameters. This may suggest different properties of GPM in the logrithmic sparsity region and the linear region.

\begin{figure}[h!]
    \centering

    \begin{subfigure}[b]{0.325\textwidth}
        \includegraphics[trim = 45mm 80mm 45mm 90mm, clip, width=\textwidth]{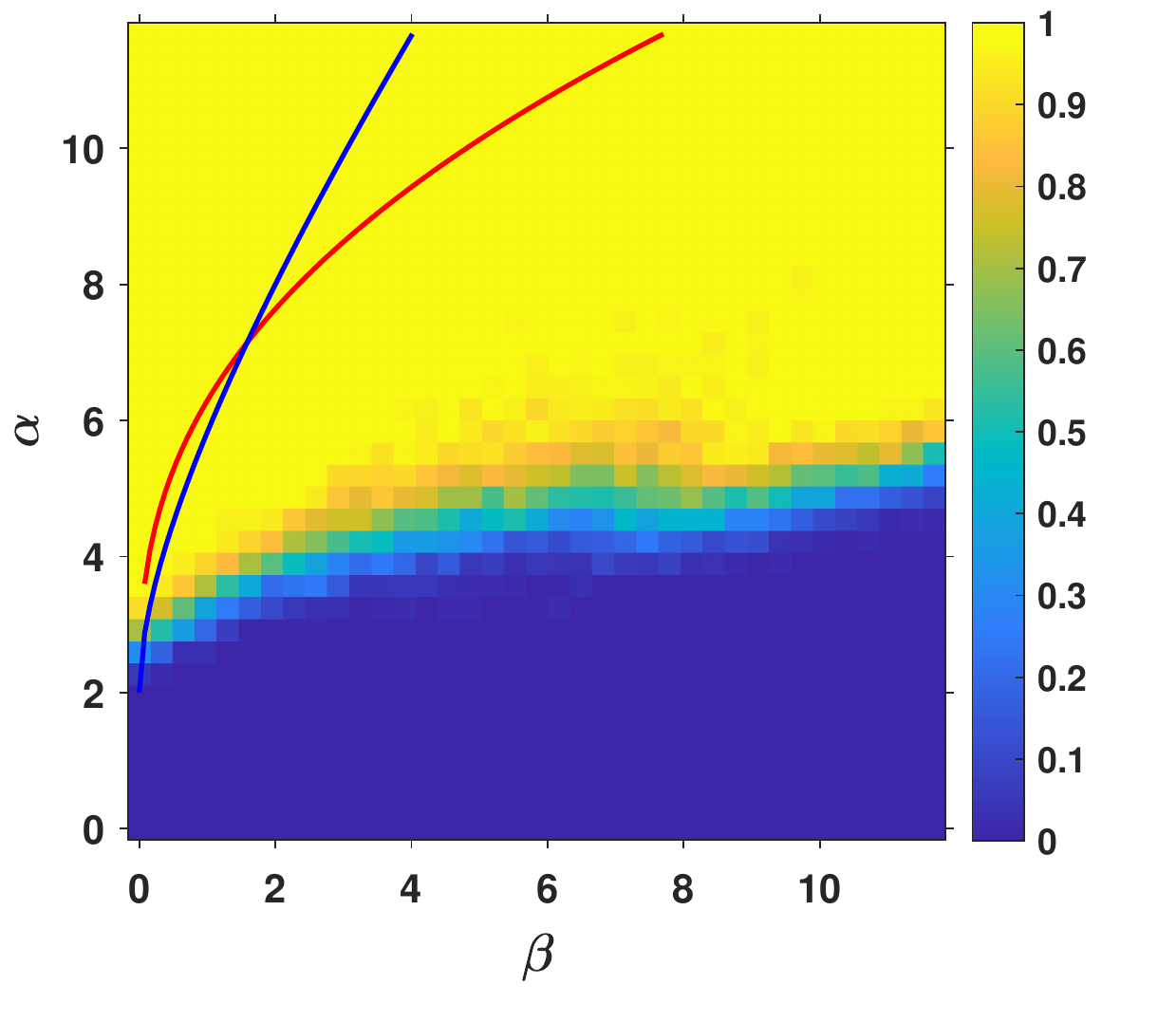}\label{fig:two_equal_a}
        \caption{GPM}
    \end{subfigure}
    \hspace{20mm}
    \begin{subfigure}[b]{0.325\textwidth}
        \includegraphics[trim = 45mm 80mm 45mm 90mm, clip, width=\textwidth]{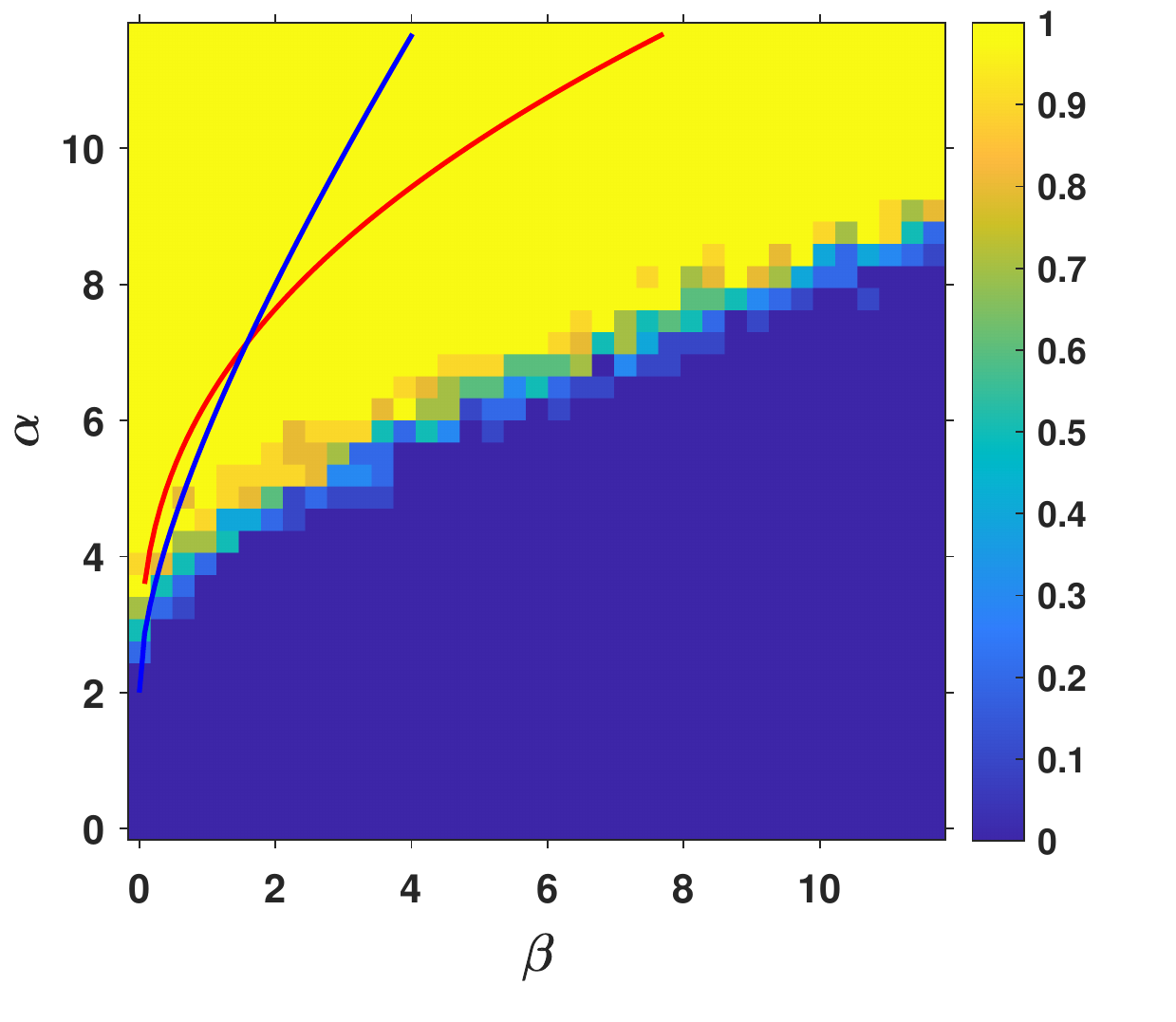}\label{fig:two_equal_b}
        \caption{SDP}
    \end{subfigure}

    \caption{Phase transition results on GPM and SDP with $\mathcal{G}=\SO(3)$ and $(n,K)=(50,2)$. The theoretical threshold for pure community detection $\sqrt{\alpha}-\sqrt{\beta}=\sqrt{K}$ is plotted in blue; the improved lower bound claimed in Theorem \ref{thm:cond-i-ii} is plotted in red.}
    \label{fig:phase-transition-comparison}
\end{figure}

Then, the phase transition of our GPM is compared with SDP proposed in \cite{fan2021joint}, which we implement with MATLAB CVX package and MOSEK solver. We only experiment on rotational scenarios considering the focus of \cite{fan2021joint}. According to the results shown in Figure \ref{fig:phase-transition-comparison}, the recovery performance of GPM notably outperforms SDP, but the boundary of transition is less sharp.

\subsection{Convergence performance and CPU time}
\begin{figure}[h!]
    \centering
    \begin{subfigure}[b]{0.325\textwidth}
        \includegraphics[trim = 35mm 80mm 45mm 90mm, clip, width=\textwidth]{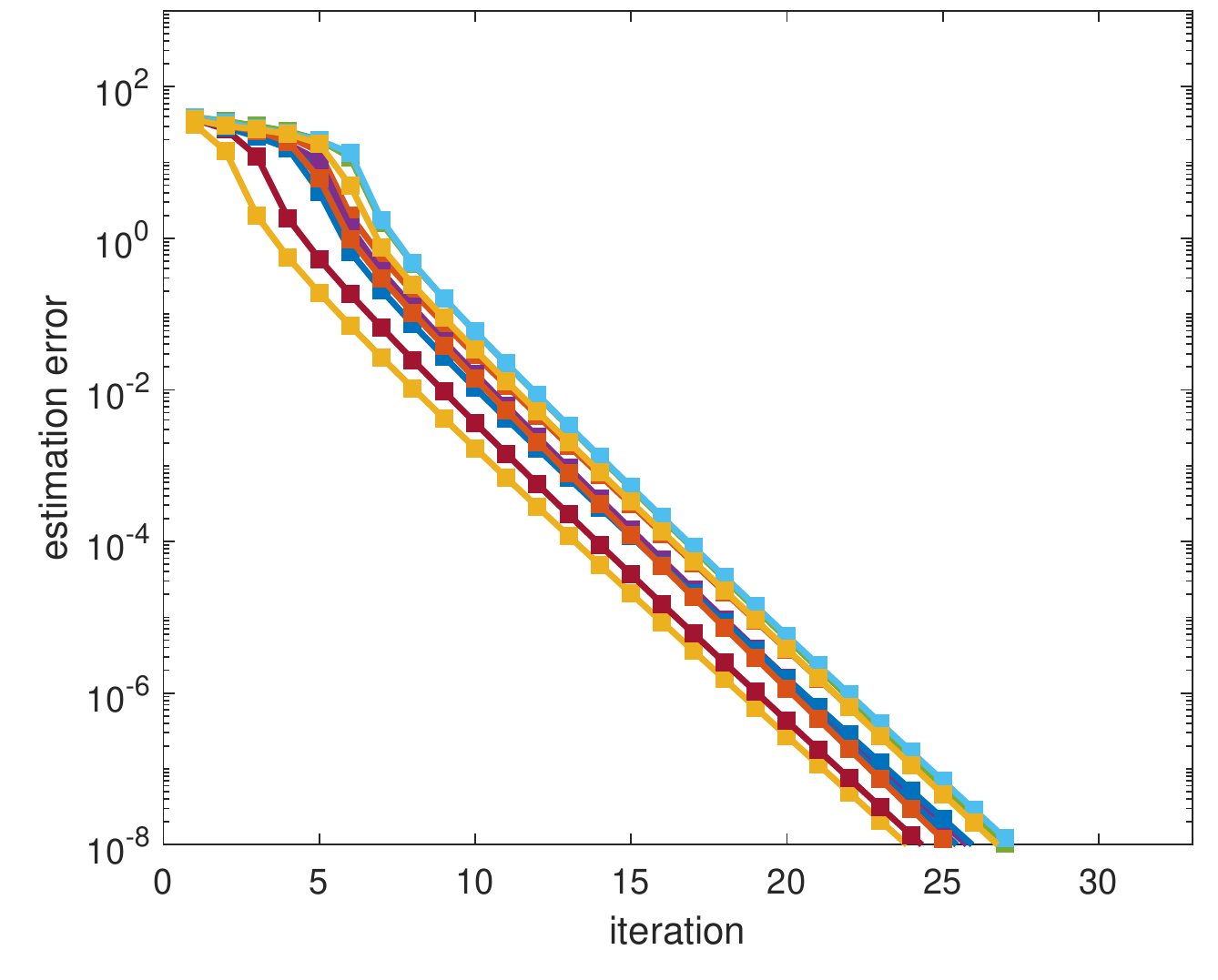}\label{fig:ps_a}
        \caption{$K=5, \alpha=15, \beta=10$}
    \end{subfigure}
    \hfill
    \begin{subfigure}[b]{0.325\textwidth}
        \includegraphics[trim = 35mm 80mm 45mm 90mm, clip, width=\textwidth]{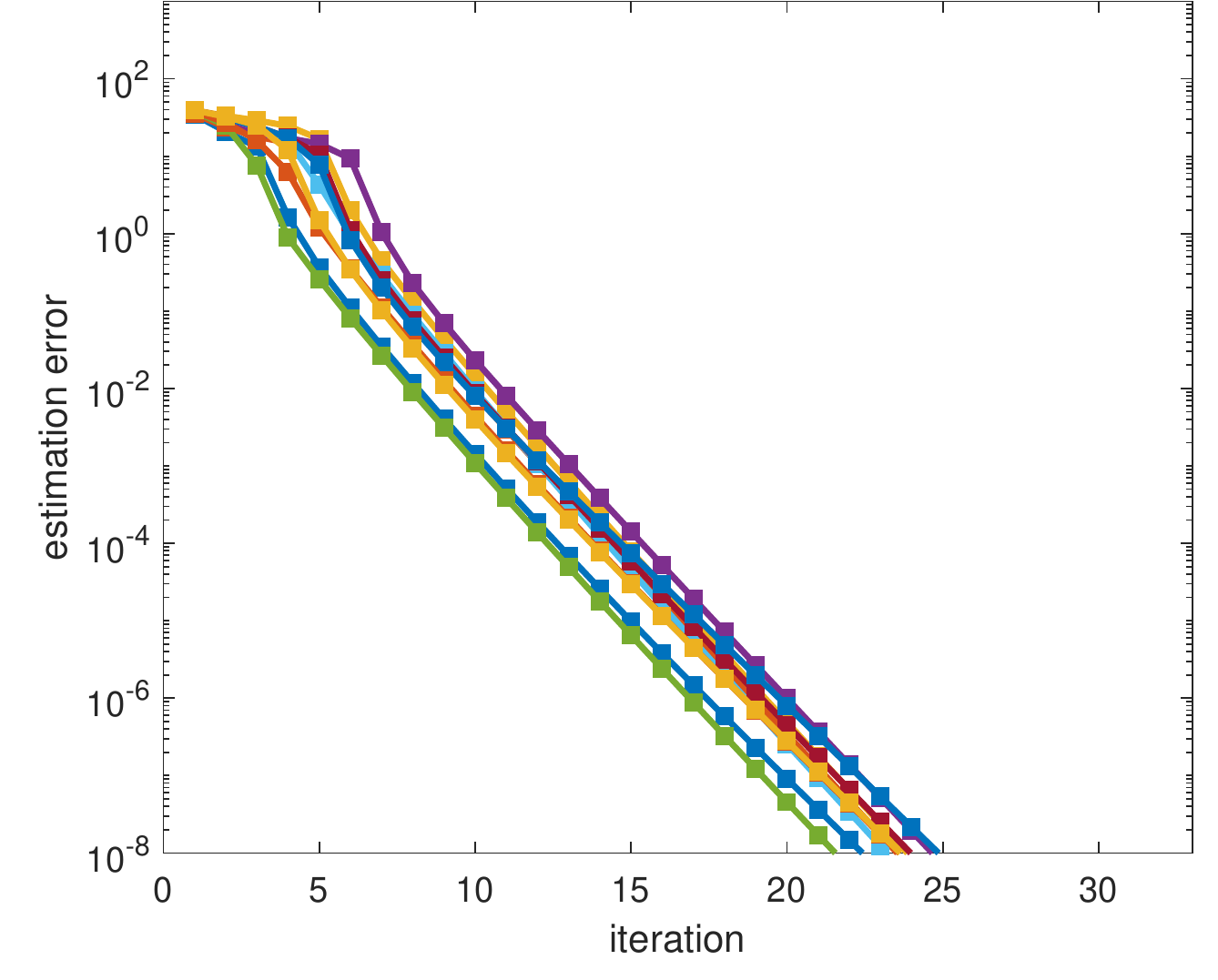}\label{fig:ps_b}
        \caption{$K=8, \alpha=25, \beta=15$}
    \end{subfigure}
    \hfill
    \begin{subfigure}[b]{0.325\textwidth}
        \includegraphics[trim = 35mm 80mm 45mm 90mm, clip, width=\textwidth]{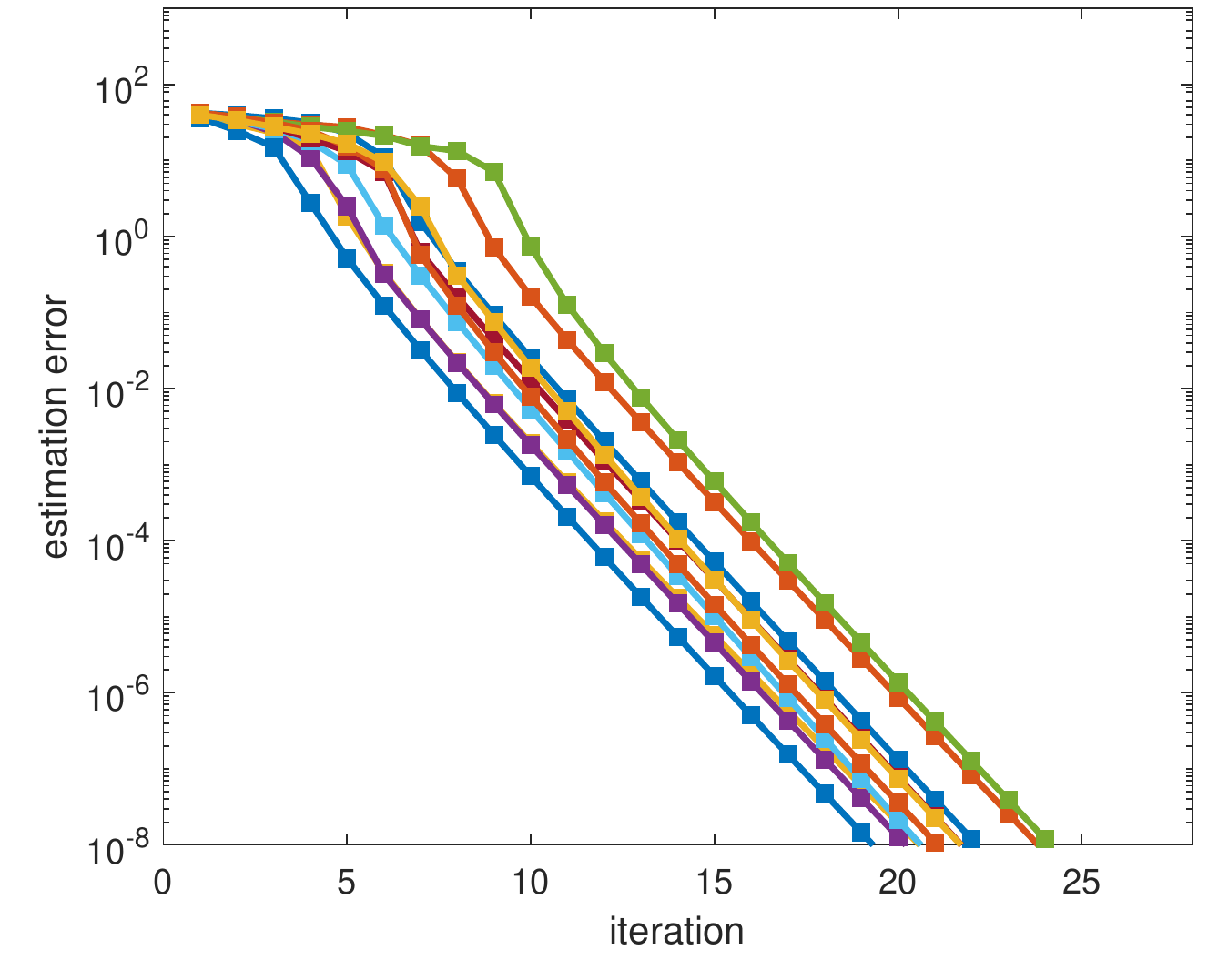}\label{fig:ps_c}
        \caption{$K=10, \alpha=35, \beta=25$}
    \end{subfigure}
    \caption{Convergence results of GPM at $\G=\OO(d)$ and $n=400$, with three different groups of parameters $(K,\alpha,\beta)=(5,15,10)$, $(8,25,15)$, and $(10,35,25)$.}
    \label{fig:convergence}
\end{figure}

The convergence performance of GPM is also studied. At $\G=\OO(d),n=400$ and three different settings of $(K,\alpha,\beta)$, we keep record of the estimation error stated in Definition \ref{def:errorV} for each iteration $\{\bm V^t\}_{t\geq 0}$ executed by GPM. Experiments are repeated for 10 times for each group of parameters to integrate the general patterns, with the resulting convergence curves plotted in Figure \ref{fig:convergence}. Typically, GPM is able to recover the ground truth within a considerably small number of iterations after a linear decay of estimation error, which again aligns with our theoretical affirmation. To further evidence the superiority of GPM in respect of its time complexity, we also test the average CPU time of both GPM and SDP on problems of different scales. All the results reported in Table \ref{tab:CPU-time} obviously demonstrate a higher time efficiency of our algorithm than SDP.

\begin{center}
\begin{table}[]
\centering
\begin{tabular}{|c|cc|ll}
\cline{1-3}
\multirow{2}{*}{parameters} & \multicolumn{2}{c|}{CPU time (s)} \\ \cline{2-3}
& \multicolumn{1}{c|}{GPM} & \multicolumn{1}{c|}{SDP} \\ 
\cline{1-3}
$n=50,\alpha=8,\beta=5$ & \multicolumn{1}{c|}{\textbf{7.37}} & \multicolumn{1}{c|}{38.69} \\
\cline{1-3}
$n=100,\alpha=15,\beta=10$  & \multicolumn{1}{c|}{\textbf{7.81}}& \multicolumn{1}{c|}{721.21} \\ 
\cline{1-3}
$n=200,\alpha=25,\beta=15$& \multicolumn{1}{c|}{\textbf{11.40}}& \multicolumn{1}{c|}{3600$+$}\\ 
\cline{1-3}
$n=400,\alpha=35,\beta=25$& \multicolumn{1}{c|}{\textbf{20.39}}& \multicolumn{1}{c|}{3600$+$}\\ 
\cline{1-3}
\end{tabular}
\caption{Comparison of the average CPU time (in seconds) between GPM and SDP at $\G=\SO(3),K=2$.}
\label{tab:CPU-time}
\end{table}
\end{center}

\section{Discussions and Conclusion}\label{sec:conc}
This paper proposes a \textit{Generalized Power Method} (GPM) to directly tackle the non-convex problem of joint community detection and group synchronization in the logarithmic sparsity region of SGBM. From the theoretical side, we give a probabilistic bound for GPM to exactly recover the ground truth in $O(n\log^2 n)$ time, as an improvement of the previous method relying on SDP \cite{fan2021joint} that typically takes $O(n^{3.5})$ time. We also propose a randomized spectral clustering method as an initializer of GPM. The corresponding analysis indicates that the initializer itself is able to yield a good approximation of the joint problem, which could also raise theoretical and practical interest.

As a generalization of the theory for SDP developed previously, our analysis of GPM applies to both orthogonal and rotational synchronization, and makes up the unresolved cases of $K\geq 3$. At the same time, our probabilistic bound breaches the information-theoretic limit for pure community detection under SBM, which implies that GPM finely exploits the additional information of group structures embedded in the joint problem. Technically, we remark that the breach originates from the concentration behavior of uniformly distributed $\SO(d)$ or $\OO(d)$ when $d\geq 2$, which is absent in the Bernoulli random variables for the pure community detection problems.

Our work also opens up further questions of the joint problem and the theoretical understanding of GPM. For example, the phase transition experiments indicate that the theoretical lower bound of parameters remains sub-optimal and calls for continuing improvements. We mainly credit this to the rather strong assumption of norm-separation on the blocks of $\bm V$, \ie, condition (\text{i}) and (\textit{ii}) in Theorem \ref{thm:main}, while GPM still performs decently without this assumption in numerical experiments. Therefore, it is hopeful to improve the lower bound by leveraging more properties of the problem and the algorithm to circumvent this assumption. Moreover, this paper takes into consideration the cases of $K=O(1)$, yet $K$ may increase with $n$ in real problems. Generalization to broader parameter regions then becomes another possible direction for future studies.

\bibliographystyle{siam}
\bibliography{refs}

\newpage
\begin{appendices}
\section{Proofs}
\subsection{Proof of Proposition \ref{prop:proj_F-runtime}}
The proof relies on the invariance of the inner product $\ip{\bm X,\bm V}$ for any $\bm X\in\R^{nd\times Kd}$ and $\bm V\in\mathcal{F}$, when transitioning the \textit{mask} $\bm H$ from $\bm V$ to $\bm X$:
\begin{lem}\label{lem:mask_transfer}
    For any $\bm X\in\R^{nd\times Kd}$, $\bm R\in\OO(d)^n$, and $\bm H\in\mathcal{H}$,
    $$\tr\Big(\bm X^\top \big(\left(\mathbf{1}_{1\times K}\otimes\bm R\right)\odot(\bm H\otimes\mathbf{1}_{d\times d})\big)\Big)=\tr\Big(\big(\bm X\odot(\bm H\otimes\mathbf{1}_{d\times d})\big)^\top \left(\mathbf{1}_{1\times K}\otimes\bm R\right)\Big).$$
\end{lem}
\begin{proof}[Proof of Proposition \ref{prop:proj_F-runtime}]
    We treat $\bm X$ as an $n$ by $K$ block matrix. By definition, 
    $$\Pi_\mathcal{F}(\bm X)=\argmin_{\bm V\in\mathcal{F}}\nm{\bm X-\bm V}_F=\argmin_{\bm V\in\mathcal{F}}\nm{\bm X}^2_F+\nm{\bm V}^2_F-\left<\bm X, \bm V\right>,$$
    where $\nm{\bm V}^2_F$ is constant because $\bm V\in\mathcal{F}$. Hence, 
    \begin{equation*}
    \Pi_\F(\bm X)=\argmax_{\bm V\in\mathcal{F}}\tr(\bm X^\top \bm V)=\argmax_{\bm R\in\OO(d)^n,\bm H\in\mathcal{H}}\tr\Big(\bm X^\top \big(\RK\odot(\bm H\otimes\mathbf{1}_{d\times d})\big)\Big).
    \end{equation*}
    Applying Lemma \ref{lem:mask_transfer},
    \begin{equation}\label{eq:proj-key}
    \begin{split}
    \Pi_\F(\bm X)&=\argmax_{\bm R\in\OO(d)^n,\bm H\in\mathcal{H}}\tr\Big(\big(\bm X\odot(\bm H\otimes\mathbf{1}_{d\times d})\big)^\top \RK)\Big)\\
    &=\argmax_{\bm R\in\OO(d)^n,\bm H\in\mathcal{H}}\tr\left(\sum_{i=1}^n(\bm X_{ie_i})^\top \bm R_i\right),
    \end{split}
    \end{equation}
    where $e_i$ is such that $\bm H_{ie_i}=1$. For any fixed $\bm H$ and every $i\in[n]$, we denote $\bm U\bm\Sigma \bm V^\top $ the SVD of $\bm X_{ie_i}$. By Proposition \ref{prop:proj-Od}, 
    \begin{equation*}
    \argmax_{\bm R_i\in\OO(d)}\tr\left(\left(\bm X_{ie_i}\right)^\top \bm R_i\right)=\bm U\bm V^\top 
    \end{equation*}
    and accordingly
    \begin{equation*}
    \max_{\bm R_i\in\OO(d)}\tr\left((\bm X_{ie_i})^\top \bm R_i\right)=\tr(\bm\Sigma)=: \bm M_{ie_i}.
    \end{equation*}
    Therefore, in order to maximize the expression in \eqref{eq:proj-key}, it suffices to find
    \begin{equation*}
        \argmax_{\bm H\in\mathcal{H}}\sum \bm M_{ie_i}=\argmax_{\bm H\in\mathcal{H}}\sum \bm M_{ij}\bm H_{ij}=\argmax_{\bm H\in\mathcal{H}}\ip{\bm H,\bm M}=\Pi_\mathcal{H}(\bm M),
    \end{equation*}
    and then to perform $\OO(d)$ projections for the selected blocks $\bm X_{ie_i}$. This validates the algorithm. Given that $K,d=\Theta(1)$, line 3 in the algorithm takes $O(n\log n)$ time according to Proposition \ref{prop:runtime-MCAP}, and all others statements take $O(n)$ time. This completes the proof.
\end{proof}

\subsection{Proof of Theorem \ref{thm:cond-i-ii}}
Firstly, we state two tail bounds for Bernoulli random variables and the random sum of uniformly distributed orthogonal matrices, respectively.
\begin{lem}[\cite{hajek2016achieving}, Lemma 2]\label{lem:Bern-tail-bound}
    Let $X\sim\mathrm{Binom}(m,\alpha\log n/n)$ for $m\in\mathbb{N},\alpha=O(1)$, where $m=\frac{n}{K}$ for some $K>0$. Let $\tau\in(0,\alpha]$. Then for a sufficiently large $n$,
    $$
    \Pr\left(X\leq\frac{\tau}{K}\log n\right)=n^{-\frac{1}{K}\left(\alpha-\tau\log\left(\frac{e\alpha}{\tau}\right)+o(1)\right)}.
    $$
\end{lem}

\begin{lem}[\cite{fan2021joint}, Theorem A.3]\label{lem:concentrate-Unif-Od}
    Suppose that $\{u_i\}_{i=1}^m$ and $\{\bm R_i\}_{i=1}^m$ are two finite random sequences independently and identically sampled from two independent distributions $\mathrm{Bern}(q)$ and $\mathrm{Unif}(\OO(d))$, respectively. Let $\bm S=\sum_{i=1}^m u_i\bm R_i$. Then, with probability at least $1-n^{-c}$,
    $$
    \nm{\bm S}\leq\sqrt{2qm(c\log n+\log 2d)}\left(\sqrt{1+\frac{c\log n+\log 2d}{18qm}}+\sqrt{\frac{c\log n+\log 2d}{18qm}}\right).
    $$
\end{lem}

\begin{rmk}
    Taking $n=Km$ and $q=\beta\log n/n$ into Lemma \ref{lem:concentrate-Unif-Od}, one can show that
    \begin{equation}\label{eq:bound-norm-S}
        \nm{\bm S}\leq\sqrt{\frac{2c\beta}{K}}\left(\sqrt{1+\frac{cK}{18\beta}}+\sqrt{\frac{cK}{18\beta}}\right)\log n
    \end{equation}
    with probability at least $1-n^{-c}$. Considering $18\beta\gg cK$, \eqref{eq:bound-norm-S} is simplified to $\nm{\bm S}\leq\sqrt{\frac{2c\beta}{K}}\log n$. This simplification is always conducted throughout the following contents.
\end{rmk}
Now we present a straightforward result on the model parameters.
\begin{lem}\label{lem:exist-tilde-tau-c}
    Suppose that the positive constants $\alpha,\beta,K, d$ are given. let $f(\tau)=\alpha-\tau\log\frac{e\alpha}{\tau}$ defined on $(0,\alpha]$. If
    \begin{numcases}{}
        &$\tauBernBound<\alpha$,\label{eq:tilde-cond1}\\
        &$\alpha-\tauBernBound\log\frac{e\alpha}{\tauBernBound}>K$,\label{eq:tilde-cond2}
    \end{numcases}
    then there exists $\tilde\tau<\alpha$, $\tilde c>1$, and $\chi>0$, such that
    \begin{numcases}{}
        &$\tilde\tau\leq\sqrt{2\tilde cK\beta}-\frac{\chi\alpha}{d}$;\label{eq:tilde-conclusion2}\\
        &$\alpha-\tilde\tau\log\frac{e\alpha}{\tilde\tau}>K$.\label{eq:tilde-conclusion3}
     \end{numcases}
\end{lem}
\begin{proof}
    One can observe that $f(\tau)=\alpha-\tau\log\left(\frac{e\alpha}{\tau}\right)$ monotonically decreases in $(0,\alpha]$. Therefore, the root $\tilde\tau$ such that $f(\tau^*)=K$ is uniquely determined in $(0,\alpha)$, and $\tau<\tau^*$ if $f(\tau)>K$ and $\tau\in(0,\alpha]$. By \eqref{eq:tilde-cond1}, there exists $c_1>1$ such that $\sqrt{2cK\beta}<\alpha$ for any $1<c\leq c_1$, and hence $\sqrt{2cK\beta}<\tau^*$. By \eqref{eq:tilde-cond2}, there exists $c_2>1$ such that 
    \begin{equation*}
        \alpha-\sqrt{2cK\beta}\log\frac{e\alpha}{\sqrt{2cK\beta}}>K
    \end{equation*}
    for $1<c\leq c_2$. Pick $\tilde c\in(1,\min\{c_1,c_2\}]$, and $\tilde\tau\in(\sqrt{2\tilde cK\beta}, \tau^*)$. \eqref{eq:tilde-conclusion2} and \eqref{eq:tilde-conclusion3} immediately follow by taking $\chi=d(\tilde\tau-\sqrt{2\tilde cK\beta})/\alpha$. 
\end{proof}

\begin{proof}[Proof for Theorem \ref{thm:cond-i-ii}]
    Denote $\bm M=\mu(\bm A\bm V^*)$. We first consider the probability of two subevents defined as follows for fixed $i,j,j'$ such that $C(i)=j\neq j'$, and then apply union bound.
    \begin{numcases}{}
        &there exists a constant $\chi>0$ such that $\bm M_{ij}-\bm M_{ij'}\geq\chi mp$;\label{eq:subevent1}\\
        &$\bm M_{ij}>\frac{\sqrt{2K\beta}}{\alpha}mp$.\label{eq:subevent2}
     \end{numcases}
    \noindent Observe that
    \begin{equation*}
        [\bm A\bm V^*]_{ij}=\sum_{k:C(k)=j}\bm A_{ik}\bm R_k^*=\sum_{\substack{k:C(k)=j\\C(k)=C(i)}}w_{ik}\bm R_i^*+\sum_{\substack{k:C(k)=j\\C(k)\neq C(i)}}u_{ik}\bm R_{ik}\bm R_k^*,
    \end{equation*}
    where $w_{ik}\sim\mathrm{Bern}(p),u_{ik}\sim\mathrm{Bern}(q), \bm R_{ik}\sim\mathrm{Unif}(\OO(d))$. In fact, the two parts in the summation are complementary, i.e.
    \begin{equation}\label{eq:AV-blocks}
        [\bm A\bm V^*]_{ij}=
        \begin{cases}
            \left(\sum_{k:C(k)=j}w_{ik}\right)\bm R_i^*, &\text{if }C(i)=j,\\
            \sum_{k:C(k)=j}u_{ik}\bm R_{ik}\bm R_k^*, &\text{otherwise.}
        \end{cases}
    \end{equation}
    Therefore, due to edge independence, $\bm M_{ij}=\nm{X\bm R_i^*}_*=dX$ where $X\sim\mathrm{Binom}(m,\alpha\log n/n)$. Likewise, denoting $\bm{S}$ the random variable as stated in Lemma \ref{lem:concentrate-Unif-Od}, $\sigma_1(\bm X_{ij'})=\nm{\bm S}$ since the distribution of $\mathrm{Unif}(\OO(d))$ is invariant under right (and left) orthogonal group actions, and consequently $\bm M_{ij'}\leq d\sigma_1(\bm X_{ij'})=d\nm{\bm S}$.
    Then both \eqref{eq:subevent1} and \eqref{eq:subevent2} are guaranteed to happen when
    \begin{numcases}{}
        &$X-\nm{\bm S}\geq \frac{\chi}{d}mp=\frac{\chi\alpha}{Kd}\log n$;\label{eq:subevents-cond1}\\
        &$\nm{\bm S}>\sqrt{\frac{2\beta}{K}}$.\label{eq:subevents-cond2}
     \end{numcases}
    \noindent With \eqref{eq:region-for-i-ii-cond1} and \eqref{eq:region-for-i-ii-cond2}, we are able to invoke Lemma \ref{lem:exist-tilde-tau-c} to find a group of parameters $\tilde\tau,\tilde c,\chi$ such that
    \begin{numcases}{}
        &$\tilde\tau<\alpha,\tilde c>1,\chi>0$;\label{eq:assump0}\\
        &$\tilde\tau\leq\sqrt{2\tilde cK\beta}-\frac{\chi\alpha}{d}$;\label{eq:assump1}\\
        &$\alpha-\tilde\tau\log\frac{e\alpha}{\tilde\tau}>K$.\label{eq:assump2}
     \end{numcases}
    \noindent Then, Lemma \ref{lem:Bern-tail-bound} indicates that
    \begin{equation}\label{eq:badevent-X}
        \Pr\left(X\geq\frac{\tilde\tau}{K}\log n\right)\geq 1-n^{-\frac{1}{K}\left(\alpha-\tilde\tau\log\frac{e\alpha}{\tilde\tau}\right)},
    \end{equation}
    while another probabilistic bound on $\nm{\bm S}$ is derived from Lemma \ref{lem:concentrate-Unif-Od}:
    \begin{equation}\label{eq:badevent-S}
        \Pr\left(\nm{\bm S}\leq \sqrt{\frac{2\tilde c\beta}{K}}\log n\right)\geq 1-n^{-\tilde c}.
    \end{equation}
    Combined with \eqref{eq:assump1}, the two events in \ref{eq:badevent-X} and \ref{eq:badevent-S} would immediately imply \eqref{eq:subevents-cond1} and \eqref{eq:subevents-cond2}, and consequently the subevents \eqref{eq:subevent1} and \eqref{eq:subevent2}. They would further establish the final proposition, given that the probability of both events stated in \eqref{eq:badevent-X} and \eqref{eq:badevent-S} is sufficiently high even after taking union bound over all $i\in[n]$ and $j'\in[K]$. However, this is guaranteed by \eqref{eq:assump0} and \eqref{eq:assump2} because, by union bound, both events hold for all $i,j'$ with probability at least
    \begin{equation*}
        1-nKn^{-\frac{1}{K}\left(\alpha-\tilde\tau\log\frac{e\alpha}{\tilde\tau}\right)}-n^{-\tilde c+1}=1-Kn^{-\frac{1}{K}\left(\alpha-\tilde\tau\log\frac{e\alpha}{\tilde\tau}-K\right)}-n^{-\tilde c+1}=1-n^{-\Omega(1)}.
    \end{equation*}
    This completes the proof. 
\end{proof}

\subsection{Proof of Lemma \ref{lem:Pi-XQ}}
\begin{proof}
    Since $\bm Q\in\mathcal{P}_K(\OO(d))$, there exists a permutation $\pi$ on $[K]$ such that $\bm Q_{\pi(i)i}\in\OO(d)$ for all $i\in[K]$, and the remaining blocks of $\bm Q$ are zero. For any $\bm X=[\bm X_{ij}]$, we have
    $$
    (\bm X\bm Q)_{ij}=\sum_{l=1}^K\bm X_{il}\bm Q_{lj}=\bm X_{i\pi(j)}\bm Q_{\pi(j)j}.
    $$
    Therefore, $\mu(\bm X\bm Q)_{ij}=\nm{\bm X_{i\pi(j)}}_*=\mu(\bm X)_{i\pi(j)}$, and $\mu(\bm X\bm Q)$ is in fact the column permutation of $\mu(\bm X)$ according to $\pi$. Denote $\bm H',\bm H$ the clustering matrices generated in the projection algorithm on the input $\bm X$ and $\bm X\bm Q$ respectively. Then, $\bm H'_{ij}=1$ if and only if $\bm H_{i\pi(j)}=1$. Now, for those $i,j$ such that $\bm H'_{ij}=1$, we have
    \begin{align}
    \Pi_\mathcal{F}(\bm X\bm Q)_{ij}&=\Pi_{\OO(d)}\left((\bm X\bm Q)_{ij}\right)=\Pi_{\OO(d)}\left(\bm X_{i\pi(j)}\bm Q_{\pi(j)j}\right)\nonumber\\
    &=\Pi_{\OO(d)}\left(\bm X_{i\pi(j)}\right)\bm Q_{\pi(j)j}=\Pi_\mathcal{F}(\bm X)_{i\pi(j)}\bm Q_{\pi(j)j}=\sum_{l=1}^K\Pi_\mathcal{F}(\bm X)_{il}\bm Q_{lj}.\nonumber
    \end{align}
    Hence $\Pi_\mathcal{F}(\bm X\bm Q)=\Pi_\mathcal{F}(\bm X)\bm Q$.
\end{proof}

\subsection{Proof of Lemma \ref{lem:fixed-point}}
\begin{proof}
    Denote $\bm H'$ and $\bm H^*$ the clustering matrices determined in the projection algorithm on the input $\bm A\bm V^*$ and $\bm V^*$, respectively. Since condition (\textit{i}) holds, $\bm H'=\bm H^*$. By \eqref{eq:AV-blocks}, $\Pi_{\OO(d)}\left((\bm A\bm V^*)_{iC(i)}\right)=\bm R_i^*$. Hence $\Pi_\mathcal{F}(\bm A\bm V^*)=\bm V^*$.
\end{proof}

\subsection{Proof of Proposition \ref{thm:Lipschitz}}
In order to establish the Lipschitz-like property of $\Pi_\mathcal{F}$, we first show that the maps $\mu$ and $\Pi_{\OO(d)}$ involved in the computation of $\Pi_\mathcal{F}$ have a similar behavior.
\begin{lem}\label{lem:Lipschitz-M}
    For any $\bm X,\bm X'\in\R^{nd\times Kd}$, $$\nm{\mu(\bm X)-\mu(\bm X')}_F\leq\sqrt{d}\nm{\bm X-\bm X'}_F.$$
\end{lem}
\begin{proof}
    For simplicity we denote $\sigma_k=\sigma_k(\bm X_{ij})$ and $\sigma'_k=\sigma_k(\bm X'_{ij})$. Then
    \begin{align}
        \left|\mu(\bm X)_{ij}-\mu(\bm X')_{ij}\right|&=\left|\sigma_1-\sigma'_1+\sigma_2-\sigma'_2+...+\sigma_d-\sigma'_d\right |\nonumber\\
        &\leq\sqrt{d}\sqrt{(\sigma_1-\sigma'_1)^2+(\sigma_2-\sigma'_2)^2+...+(\sigma_d-\sigma'_d)^2}\nonumber\\
        &\leq\sqrt{d}\nm{\bm X_{ij}-\bm X'_{ij}}_F,\nonumber
    \end{align}
    where Mirsky's inequality \cite{Stewart90perturbationtheory} yields the final step. Summing over the indices yields the desired result.
\end{proof}
\begin{lem}[\cite{liu2020unified}, Lemma 2]\label{lem:Lipschitz-SO}
    If $\bm X_{ij}=\eta\bm R$ where $\eta>0$ and $\bm R\in\OO(d)$, then 
    $$\nm{\Pi_{\OO(d)}(\bm X_{ij})-\Pi_{\OO(d)}(\bm X'_{ij})}_F\leq\frac{2}{\eta}\nm{\bm X_{ij}-\bm X'_{ij}}_F$$ for any $\bm X'_{ij}\in\R^{d\times d}$.
\end{lem}
\begin{proof}[Proof of Proposition \ref{thm:Lipschitz}]
    Let $\Pi_\F(\bm X)$ incoporate a community structure $\bm H$ and $\Pi_\F(\bm X')$ incoporate $\bm H'$. Then one can observe
    \begin{equation*}
        \nm{\Pi_\mathcal{F}(\bm X)-\Pi_\mathcal{F}(\bm X')}_F^2=d\nm{\bm H-\bm H'}_F^2+\sum_{j}\sum_{i\in\mathcal{I}_j\cap\mathcal{I}'_j}\nm{\Pi_{\OO(d)}(\bm X_{ij})-\Pi_{\OO(d)}(\bm X'_{ij})}_F^2.
    \end{equation*}
    By lemma 3 in \cite{wang2021optimal}, lemma \ref{lem:Lipschitz-M}, and lemma \ref{lem:Lipschitz-SO},
    \begin{align}
        \nm{\Pi_\mathcal{F}(\bm X)-\Pi_\mathcal{F}(\bm X')}_F^2&\leq \frac{4d}{\delta^2}\nm{\bm M-\bm M'}_F^2+\sum_{j}\sum_{i\in\mathcal{I}_j\cap\mathcal{I}'_j}\nm{\Pi_{\OO(d)}(\bm X_{ij})-\Pi_{\OO(d)}(\bm X'_{ij})}_F^2\nonumber\\
        &\leq \frac{4d^2}{\delta^2}\nm{\bm X-\bm X'}_F^2+\frac{4}{\eta^2}\sum_{j}\sum_{i\in\mathcal{I}_j\cap\mathcal{I}'_j}\nm{\bm X_{ij}-\bm X'_{ij}}_F^2\nonumber\\
        &\leq \left(\frac{4d^2}{\delta^2}+\frac{4}{\eta^2}\right)\nm{\bm X-\bm X'}_F^2.\nonumber
    \end{align}
    Taking $\delta=\chi mp$ and $\eta=\frac{\tauBernBound}{\alpha}mp$ yields the result.
\end{proof}

\subsection{Proof of Proposition \ref{prop:A-pVVT}}
\begin{proof}
    This is a direct generalization of Lemma 3.6 in \cite{fan2021joint} when $K>2$ and the constraint is relaxed from $\mathcal{SO}(d)$ to $\mathcal{O}(d)$. We apply similar notations. Observe that
    \begin{equation*}
        \bm{S}_\text{out}=\begin{pmatrix}
            \mathbf{0}&\bm{S}_{12}&\bm{S}_{13}&\cdots&\bm{S}_{1K}\\
            \bm{S}_{12}^\top&\mathbf{0}&\bm{S}_{23}&\cdots&\bm{S}_{2K}\\
            \bm{S}_{13}^\top&\bm{S}_{23}^\top&\mathbf{0}&\cdots&\bm{S}_{3K}\\
            \vdots&\vdots&\vdots&\ddots&\vdots\\
            \bm{S}_{1K}^\top&\bm{S}_{2K}^\top&\bm{S}_{3K}^\top&\cdots&\mathbf{0}
        \end{pmatrix}=\sum_{j>i}\begin{pmatrix}
            &&&&\\
            &&&\bm{S}_{ij}&\\
            &&&&\\
            &\bm{S}_{ij}^\top&&&\\
            &&&&
        \end{pmatrix}.
    \end{equation*}
    Therefore, $\nm{\bm{S}_\text{out}}\leq\sum_{j>i}\nm{\bm{S}_{ij}}$, and likewise $\nm{\bm{S}_\text{in}}\leq\sum_{i}\nm{\bm{S}_{ii}}$. Following the argument therein, the result can be established by union bound.
\end{proof}

\subsection{Proof of Proposition \ref{prop:Z-I-bound}}
\begin{proof}
    We denote for simplicity $\bm V^e=\bm V^*\bm Q^t$, and $\bm Z=\frac{1}{m}\bm V^{e\top}\bm V^t$. Recall that $\bm V^e$ is the optimal approximation of $\bm V^t$, so any per-cluster orthogonal transformation never yields a smaller difference. Specifically,
    \begin{equation*}
        \nm{\bm V^t-\bm V^e}_F^2=\min_{\bm P}\nm{\bm V^t-\bm V^e\bm P^\top }_F^2
    \end{equation*}
    subject to
    \begin{equation*}
        \bm P=\begin{pmatrix}
            \bm P_1&&&\\
            &\bm P_2&&\\
            &&\ddots&\\
            &&&\bm P_K
        \end{pmatrix}\in\mathcal{P}_K(\OO(d)),\bm P_i\in\OO(d).
    \end{equation*}
    However, note that
    \begin{align}
        \min_{\bm P}\nm{\bm V^t-\bm V^e\bm P^\top }_F^2&=\sum_{i=1}^K\min_{\bm P_i\in\OO(d)}\left(\sum_{j\in\mathcal{I}^e_i\cap\mathcal{I}_i}\nm{\bm V^t_{ji}-\bm V_{ji}^e\bm P_i^\top }_F^2+\sum_{j\in\mathcal{I}^e_i\cup\mathcal{I}_i/\mathcal{I}^e_i\cap\mathcal{I}_i}\nm{\bm V^t_{ji}-\bm V_{ji}^e\bm P_i^\top }_F^2\right)\nonumber\\
        &=\sum_{i=1}^K\left(2Kmd-2\max_{\bm P_i\in\OO(d)}\sum_{j\in\mathcal{I}^e_i\cap \mathcal{I}_i}\ip{\bm V^t_{ji}, \bm V_{ji}^e\bm P_i^\top }\right)\nonumber\\
        &=\sum_{i=1}^K\left(2Kmd-2\max_{\bm P_i\in\OO(d)}\ip{\bm V^t_{ji}, \sum_{j\in\mathcal{I}^e_i\cap \mathcal{I}_i}\bm V_{ji}^e\bm P_i^\top }\right)\nonumber\\
        &=\sum_{i=1}^K\left(2Kmd-2m\max_{\bm P_i\in\OO(d)}\ip{\bm P_i, \bm Z_{ii}^\top }\right),\nonumber
    \end{align}
    we obtain $\bm P_i=\Pi_{\OO(d)}\left(\bm Z_{ii}^\top \right)$ (and $\bm P_i=\bm I_d$ at the same time). By Proposition \ref{prop:proj-Od}, we have
    \begin{align}
    &\tr(\bm Z_{ii})=\ip{\bm I,\bm Z_{ii}^\top }=\ip{\Pi_{\OO(d)}\left(\bm Z_{ii}^\top \right), \bm Z_{ii}^\top }=\sum_{k=1}^d\sigma_k\left(\bm Z_{ii}\right);\label{eq:tr-sum-of-sigma}\\
    &\nm{\bm Z_{ii}}_F^2=\nm{\Pi_{\OO(d)}\left(\bm Z_{ii}^\top \right)\bm Z_{ii}}_F^2=\sum_{k=1}^d\sigma_k\left(\bm Z_{ii}\right)^2;\label{eq:norm-square-sum-of-sigma}\\
    &\sigma_k(\bm Z_{ii})\in[0,1].\label{eq:sigma-range}
    \end{align}
    
    \noindent We now claim that
    \begin{equation*}
        \nm{\bm Z-\bm I_{Kd}}_F^2\leq\left(1+\frac{1}{d}\right)[Kd-\tr(\bm Z)]^2.
    \end{equation*}
    To this end, observe that
    \begin{align}
        \nm{\bm Z-\bm I_{Kd}}_F^2&=\sum_{i=1}^K\left(\nm{\bm Z_{ii}}_F^2+\sum_{j\neq i}\nm{\bm Z_{ij}}_F^2-2\tr(\bm Z_{ii})+d\right)\nonumber\\
        &\leq \sum_{i=1}^K\left(\nm{\bm Z_{ii}}_F^2+\left(\sqrt{d}-\nm{\bm Z_{ii}}_F\right)\sum_{j\neq i}\nm{\bm Z_{ij}}_F-2\tr(\bm Z_{ii})+d\right)\nonumber\\
        &\leq \sum_{i=1}^K\left(\nm{\bm Z_{ii}}_F^2+\left(\sqrt{d}-\nm{\bm Z_{ii}}_F\right)^2-2\tr(\bm Z_{ii})+d\right).\label{eq:Z-Inorm-bound}
    \end{align}
    By \eqref{eq:tr-sum-of-sigma}, \eqref{eq:norm-square-sum-of-sigma}, and \eqref{eq:sigma-range}, we have
    \begin{align}
        \nm{\bm Z_{ii}}_F^2-2\tr(\bm Z_{ii})+d&=\sigma_1^2+...+\sigma_d^2-2(\sigma_1+...+\sigma_d)+d\nonumber\\
        &\leq\left((1-\sigma_1)+...+(1-\sigma_d)\right)^2=\left(d-\tr(\bm Z_{ii})\right)^2,\label{eq:Zii-bound}
    \end{align}
    and
    \begin{equation}\label{eq:Zij-bound}
        \left(\sqrt{d}-\nm{\bm Z_{ii}}_F\right)^2\leq\left(\sqrt{d}-\frac{\tr(\bm Z_{ii})}{\sqrt{d}}\right)^2=\frac{1}{d}\left(d-\tr(\bm Z_{ii})\right)^2.
    \end{equation}
    Summing \eqref{eq:Zii-bound} and \eqref{eq:Zij-bound} over $i$, \eqref{eq:Z-Inorm-bound} yields
    \begin{align}
        \nm{\bm Z-\bm I_{Kd}}_F^2&\leq\left(1+\frac{1}{d}\right)\sum_{i=1}^K(d-\tr(\bm Z_{ii}))^2\leq\left(1+\frac{1}{d}\right)\left(\sum_{i=1}^Kd-\tr(\bm Z_{ii})\right)^2\nonumber\\
        &=\left(1+\frac{1}{d}\right)\left(Kd-\tr(\bm Z)\right)^2,\nonumber
    \end{align}
    which validates our claim. Finally, since $\nm{\bm V^t-\bm V^e}_F^2=2mKd-2m\tr(\bm Z)$, we have
    \begin{equation*}
    \begin{split}
        m\nm{\bm Z-\bm I_{Kd}}_F\leq\frac{1}{2}\sqrt{1+\frac{1}{d}}\nm{\bm V^t-\bm V^e}_F^2\leq\frac{1}{2}\sqrt{1+\frac{1}{d}}\frac{\sqrt{m}}{\rho}\nm{\bm V^t-\bm V^e}_F
    \end{split}
    \end{equation*}
    for $\rho>0$.
\end{proof}

\subsection{Proof of Proposition \ref{prop:tightness}}
\begin{proof}
    We first observe that the community structures of $\bm V$ and $\bm V^*$ are identical up to some permutation when $\epsilon_{\OO(d)}(\bm V)<\sqrt{2}$. Otherwise, at least one node falls in a erroneous cluster and $\epsilon_{\OO(d)}(\bm V)\geq\sqrt{2d}\geq\sqrt{2}$. Now, without loss of generality, we may identify the community structure of $\bm V^*$ with that of $\bm V$. Then no permutation is required to present the equivalence class of $\bm V^*$, hence
    \begin{align}
        \epsilon_{\OO(d)}(\bm V)&=\min_{\bm W\in\mathrm{bdiag}\left(\OO(d)^K\right)}\nm{\bm V-\bm V^*\bm W}_F,\label{eq:err-Od}\\
        \epsilon_{\SO(d)}(\mathcal{R}(\bm V))&=\min_{\bm W\in\mathrm{bdiag}\left(\SO(d)^K\right)}\nm{\mathcal{R}(\bm V)-\bm V^*\bm W}_F.\label{eq:err-SOd}
    \end{align}

    Our second observation is that no two group elements in the same cluster of $\bm V$, say $\bm R_i$ and $\bm R_j$, belong to $\SO(d)$ and $\SO^-(d)$ respectively. To see this, consider $\bm V^e:=\bm V^*\bm W$ where $\bm W\in\mathrm{bdiag}(\OO(d)^K)$ is arbitrary. Observe that no two group elements in the same cluster of $\bm V^e$ belong to $\SO(d)$ and $\SO^-(d)$ respectively, because $\bm W$ exerts a unified group action on each cluster of $\bm V^*\in\mathcal{E}$. If the same does not hold for $\bm V$, there must exist some $i\in[n]$ such that $\bm R_i\in\SO(d),\bm R^e_i\in\SO^-(d)$, or $\bm R_i\in\SO^-(d),\bm R^e_i\in\SO(d)$. In both cases, $\nm{\bm R_i-\bm R^e_i}_F\geq\sqrt{2}$, $\forall\bm W\in\mathrm{bdiag}(\OO(d)^K)$, hence 
    $$
    \epsilon_{\OO(d)}(\bm V)=\min\nm{}_F\geq\sqrt{2}.
    $$
    Therefore, when $\epsilon_{\OO(d)}(\bm V)<\sqrt{2}$, the rounding procedure $\mathcal{R}(\bm V)=\bm V\bm T$, where
    \begin{equation*}
        \bm T=\begin{pmatrix}
            D_1\bm I_d&&&\\
            &D_2\bm I_d&&\\
            &&\ddots&\\
            &&&D_K\bm I_d
        \end{pmatrix}\in\mathrm{bdiag}(\OO(d)^K),\;D_k=\det(\bm R_i)\;\forall i\in\mathcal{I}_k.
    \end{equation*}
    This together with \eqref{eq:err-SOd} gives
    \begin{equation}\label{eq:err-SOd-sqrt2}
        \epsilon_{\SO(d)}(\mathcal{R}(\bm V))=\min_{\bm W\in\mathrm{bdiag}\left(\SO(d)^K\right)}\nm{\bm V\bm T-\bm V^*\bm W}_F.
    \end{equation}
    Moreover, since $\bm T\in\mathrm{bdiag}(\OO(d)^n)$, \eqref{eq:err-Od} gives
    \begin{equation}\label{eq:err-Od-sqrt2}
        \epsilon_{\OO(d)}(\bm V)=\min_{\bm W\in\mathrm{bdiag}(\OO(d)^K)}\nm{\bm V-\bm V^*\bm W}_F=\min_{\bm W\in\mathrm{bdiag}(\OO(d)^K)}\nm{\bm V\bm T-\bm V^*\bm W}_F.
    \end{equation}
    Denote $\bm W^*$ a minimizer of \eqref{eq:err-Od-sqrt2}, \ie,
    \begin{equation*}
        \bm W^*=\argmin_{\bm W\in\mathrm{bdiag}(\OO(d)^K)}\nm{\bm V\bm T-\bm V^*\bm W}_F.
    \end{equation*}
    Then $\nm{\bm V\bm T-\bm V^*\bm W^*}_F<\sqrt{2}$. Since $\bm V\bm T\in\mathcal{E}$, it follows from an argument similar to the second observation that $\bm W^*\in\mathrm{bdiag}(\SO(d)^K)$, lest the estimation error exceeds $\sqrt{2}$. Therefore, $\bm W^*$ also minimizes \eqref{eq:err-SOd-sqrt2}. We then establish the equality $\epsilon_{\SO(d)}(\bm V\bm T)=\epsilon_{\OO(d)}(\bm V)$.
\end{proof}

\subsection{Proof of Proposition \ref{prop:rand-spec-init}}
We make use of the following variant of Davis-Kahan theorem on the distance between eigenspaces of two real symmetric matrices.
\begin{prop}[Davis-Kahan, \cite{yu2014DKtheorem}]\label{lem:DK}
    Let $\bm M, \bm M^*\in\R^{N\times N}$ be symmetric matrices with eigenvalues $\lambda_1\geq...\geq\lambda_N$ and $\lambda^*_1\geq...\geq\lambda^*_N$, respectively. For any integers $k,l$ such that $1\leq k\leq l\leq N$, let $\bm U=\mathrm{eigs}_{[k:l]}(\bm M),\bm U^*=\mathrm{eigs}_{[k:l]}(\bm M^*)$. Suppose that $\min\{\lambda^*_{k-1}-\lambda^*_k,\lambda^*_l-\lambda^*_{l+1}>0\}$, where $\lambda_0=+\infty,\lambda_{N+1}=-\infty$. Then, there exists $\bm Q^*\in\OO(l-k+1)$ such that
    \begin{equation*}
        \nm{\bm U-\bm U^*\bm Q^*}_F\leq\frac{2\sqrt{2}\sqrt{l-k+1}\nm{\bm M-\bm M^*}}{\min\{\lambda^*_{k-1}-\lambda^*_k,\lambda^*_l-\lambda^*_{l+1}\}}.
    \end{equation*}
\end{prop}

\begin{proof}[Proof of Theorem \ref{thm:cond-iii}]We prove the existence of such an algorithm by showing that Algorithm \ref{alg:spec-init} does satisfy all the desired properties. Observe that $\bm V^*/\sqrt{m}$ are the leading eigenvectors of $p\bm V^*\bm V^{*\top}$ with eigenvalues $pm$, while the other eigenvalues of $p\bm V^*\bm V^{*\top}$ are all zero. By Lemma \ref{lem:DK}, there exists $\bm Q^*\in\mathcal{O}(Kd)$ such that
\begin{equation*}
    \nm{\hU-\frac{1}{\sqrt{m}}\bm V^*\bm Q^*}_F\leq\frac{2\sqrt{2Kd}}{pm}\nm{\bm A-p\bm V^*\bm V^{*\top}}.
\end{equation*}
Denote $\bm\Phi=\frac{1}{\sqrt{m}}\bm V^*\bm Q^*$. By Proposition \ref{prop:A-pVVT}, for sufficiently large $n$, there exists $c_4>0$ such that
\begin{equation}\label{eq:U-Phi-norm-bound}
    \nm{\hU-\bm\Phi}_F\leq\frac{2\sqrt{2Kd}}{pm}\left(c_1\sqrt{qm}+c_2\sqrt{pm}+c_3\sqrt{\log n}\right)<\frac{c_4}{\sqrt{\log m}}.
\end{equation}
Also, by direct calculation,
\begin{equation*}
    \bm\Phi_{v\times}=\frac{1}{\sqrt{m}}\bm R^*_{v}\bm Q^*_{C^*(v)\times},
\end{equation*}
and it is a direct consequence that for all $v\in[n]$,
\begin{equation}\label{eq:norm-Phi-v}
    \nm{\bm\Phi_{v\times}}_F=\frac{d}{m}.
\end{equation}
Moreover, for $v, u$ belonging to the same ground truth cluster, we have
\begin{equation*}
    \bm\Phi_{v\times}\bm\Phi_{u\times}^\top =\frac{1}{m}\bm R_v^*\bm R_u^{*\top}.
\end{equation*}
Lemma \ref{lem:Lipschitz-SO} then implies
\begin{equation}\label{eq:spec-init-O}
    \nm{\Pi_{\mathcal{O}(d)}\left(\hU_{v\times}\hU_{u\times}^\top \right)-\bm R_v^*\bm R_u^{*\top}}_F\leq 2m\nm{\hU_{v\times}\hU_{u\times}^\top -\bm\Phi_{v\times}\bm\Phi_{u\times}^\top }_F.
\end{equation}
Now we consider $u=\tau_i$ and $v\in\mathcal{I}^0_i\cap\mathcal{I}^*_{\pi(i)}$. Suppose the following conditions hold for all $i\in[n]$, whose validity with high probability will be proved at the end of this section:
\begin{numcases}{}
    &$\tau_i\in\mathcal{I}^*_{\pi(i)}$;\label{eq:tau-i-assumption1}\\
    &there exists a constant $c_5>0$ such that $\nm{\hU_{\tau_i\times}-\bm\Phi_{\tau_i\times}}_F\leq\sqrt{\frac{1}{(8Kd+c_5)\rho^2 m}}$.\label{eq:tau-i-assumption2}
\end{numcases}
Then \eqref{eq:spec-init-O} yields
\begin{equation*}
    \nm{\bm R^0_v-\bm R^*_v\bm R^{*\top}_{\tau_i}}_F\leq 2m\nm{\hU_{v\times}\hU_{\tau_i\times}^\top -\bm\Phi_{v\times}\bm\Phi_{\tau_i\times}^\top }_F.
\end{equation*}
Therefore, if we denote $\tau_{C^*(v)}=\delta(v)$,
\begin{align}
    &\sum_{i=1}^K\sum_{v\in\mathcal{I}_i^0\cap \mathcal{I}^*_{\pi(i)}}\nm{\bm R^0_v-\bm R^*_v\bm R^{*\top}_{\tau_i}}_F^2\leq 4m^2\sum_{v=1}^n\nm{\hU_{v\times}\hU_{\delta(v)\times}^\top -\bm\Phi_{v\times}\bm\Phi_{\delta(v)\times}^\top }_F^2\nonumber\\
    \leq &8m^2\left( \sum_{v=1}^n\nm{\left(\hU_{v\times}-\bm\Phi_{v\times}\right)\hU_{\delta(v)\times}^\top }_F^2 + \sum_{v=1}^n\nm{\left(\hU_{\delta(v)\times}-\bm\Phi_{\delta(v)\times}\right)\bm\Phi_{v\times}^\top }_F^2 \right)\nonumber\\
    \leq & 8m^2\max_{v\in[n]}\nm{\hU_{\delta(v)\times}^\top}_F^2\times\sum_{v=1}^n\nm{\hU_{v\times}-\bm\Phi_{v\times}}_F^2+8m^2\frac{d}{m}\sum_{v=1}^n\nm{\hU_{\delta(v)\times}-\bm\Phi_{\delta(v)\times}}_F^2,\nonumber
\end{align}
where the second inequality follows from triangle inequality, and the third from \eqref{eq:norm-Phi-v} and the inequality $\nm{\bm X\bm Y}_F\leq\nm{\bm X}_F\nm{\bm Y}_F$. Apply triangle inequality to \eqref{eq:tau-i-assumption2}, we have $\nm{\hU_{\delta(v)\times}}_F^2\leq\frac{c_6}{m}$ for some constant $c_6>0$. \eqref{eq:tau-i-assumption2} again implies
\begin{equation*}
    \begin{split}
        \sum_{i=1}^K\sum_{v\in\mathcal{I}_i^0\cap \mathcal{I}^*_{\pi(i)}}\nm{\bm R^0_v-\bm R^*_v\bm R^{*\top}_{\tau_i}}_F^2
        \leq \frac{8c_4^2c_6m}{\log m}+8m^2Kd\frac{1}{(8Kd+c_5)\rho^2 m}<\frac{c_7m}{\rho^2},
    \end{split}
\end{equation*}
where $0<c_7<1$ is a constant. This yields
\begin{align}
\epsilon(\bm V^0)^2&\leq\frac{Cnd}{\log n}+\sum_{i=1}^K\sum_{v\in\mathcal{I}_i^0\cap \mathcal{I}^*_{\pi(i)}}\nm{\bm R^0_v-\bm R^*_v\bm R^{*\top}_{\tau_i}}_F^2\nonumber\\
&\leq \frac{Cnd}{\log n}+\frac{c_7m}{\rho^2}<\frac{m}{\rho^2},\nonumber
\end{align}
which completes the proof.

We now show that \eqref{eq:tau-i-assumption1} and \eqref{eq:tau-i-assumption2} simultaneously hold with probability at least $1-(\log n)^{-\Omega(1)}$. Firstly, according to \eqref{eq:H0-err}, there exists a permutation $\pi$ of the set $[K]$ such that 
$$\left|\mathcal{I}_i^0\cap\mathcal{I}^*_{\pi(i)}\right|\geq m-\frac{CKm}{2\log Km}$$
for arbitrary $i\in[K]$. Therefore, for any fixed $i$, $\tau_i$ picked in algorithm \ref{alg:spec-init} satisfy \eqref{eq:tau-i-assumption1} with probability at least
\begin{equation}\label{eq:tau-i-prob-1}
    1-\frac{CK}{2\log Km}.
\end{equation}
Secondly, for any size-$m$ set $T\subset [n]$, the size of the subset 
$$\left\{t\in T:\nm{\hU_{t\times}-\bm\Phi_{t\times}}_F^2\leq\frac{1}{(8Kd+c_5)\rho^2m}\right\}$$
is at least $m-\frac{m}{\sqrt{\log m}}$. Otherwise \eqref{eq:U-Phi-norm-bound} is contradicted since
\begin{equation*}
    \nm{\hU-\bm\Phi}_F^2\geq\sum_{t\in T}\nm{\hU_{t\times}-\bm\Phi_{t\times}}_F^2>\frac{m}{\sqrt{\log m}}\frac{1}{(8Kd+c_5)\rho^2m}>\frac{c_4^2}{\log m}
\end{equation*}
for a sufficiently large $m$. Therefore, for any fixed $i$, \eqref{eq:tau-i-assumption2} holds with probability at least 
\begin{equation}\label{eq:tau-i-prob-2}
    1-\frac{1}{\sqrt{\log m}}.
\end{equation}
By \eqref{eq:tau-i-prob-1}, \eqref{eq:tau-i-prob-2} and union bound, \eqref{eq:tau-i-assumption1} and \eqref{eq:tau-i-assumption2} simultaneously hold for all $i$ with probability at least
$$1-K\left(\frac{CK}{2\log Km}+\frac{1}{\sqrt{\log m}}\right)=1-(\log n)^{-\Omega(1)}.$$
\end{proof}

\end{appendices}

\end{document}